\documentclass[11pt]{article}

\usepackage{fullpage}
\usepackage{amssymb,amsthm}
\usepackage{amsmath}
\usepackage{amsfonts}
\usepackage{latexsym}
\usepackage{amsxtra}
\usepackage{graphicx}
\usepackage{caption}
\usepackage{subcaption}
\usepackage{wrapfig}
\usepackage{mathrsfs}
\usepackage{epsfig}
\usepackage{times}
\usepackage{makeidx}
\usepackage[titles]{tocloft}
\usepackage{cite}

\newcommand{\heikodetail}[1]{}
\renewcommand{\proof}{\noindent{\sc Proof:}\hskip 1.1em}
\renewcommand{\qed}{\hfill\mbox{$\Box$}\\}
\def\yy{\mbox{$\spadesuit$}}
\newcommand{\invisible}[1]{\par ($\spadesuit$ \emph{hidden comments in \TeX{}
file\/})}

\newcommand{\omitted}[1]{}

\newcommand{\cC}{{\mathcal{C}}}
\newcommand{\U}{{\mathscr{U}}}

\newtheorem{lemma}{\bf Lemma}[section]
\newtheorem{theorem}[lemma]{\bf Theorem}
\newtheorem{proposition}[lemma]{\bf Proposition}
\newtheorem{corollary}[lemma]{\bf Corollary}

\theoremstyle{definition}
\newtheorem{definition}[lemma]{\bf Definition}
\newtheorem{REMARK}[lemma]{\bf Remark}

\newcommand{\Fo}{\,\,\,\text{for }\,\,}
\newcommand{\Foa}{\,\,\,\text{for all }\,\,}

\newcommand{\AND}{\,\,\,\text{and }\,\,}

\newcommand\Reals{{\mathbb R}}
\newcommand\R{{\mathbb R}}
\newcommand\Z{{\mathbb Z}}

\newcommand\N{{\mathbb N}}

\newcommand{\bbbr}{\Reals}

\newcommand\dist{\mathop{\rm dist}\nolimits}
\newcommand\diam{\mathop{\rm diam}\nolimits}

\newcommand\ang{\mathop{\mbox{$<\!\!\!)$}}\nolimits}

\newcommand{\xx}{\mbox{$\clubsuit$}}
\newcommand{\F}{\mathscr{F}}
\newcommand{\E}{\mathscr{E}}

\newcommand{\M}{\mathscr{M}}
\newcommand{\I}{\mathscr{I}}

\newcommand{\eps}{\epsilon}
\renewcommand{\H}{\mathscr{H}}

\def\emob{\E_{\text{M\"{o}b}}}
\def\esym{\E^{\text{sym}}}
\def\acn{\text{acn}}

\def\circle{\text{circle}}

\begin{document}

\title{On some knot energies involving Menger curvature}

\author{Pawe\l{} Strzelecki, Marta Szuma\'{n}ska, Heiko von der Mosel}


\maketitle

\frenchspacing

\begin{abstract}

We investigate knot-theoretic properties of geometrically defined curvature energies such as integral Menger curvature.
Elementary radii-functions, such as the circumradius of three points, generate a~family of knot energies guaranteeing
self-avoidance and a~varying degree of higher regularity of  finite energy curves. All of these energies turn out to be
charge, minimizable in given isotopy classes, tight and strong. Almost all distinguish between knots and unknots, and some
of them can be shown to be uniquely minimized by round circles. Bounds on the stick number and the average crossing number,
some non-trivial global lower bounds, and unique minimization by circles upon compaction complete the picture.
\vspace{3mm}

\centering{Mathematics Subject Classification  (2010): 49Q10, 53A04,   57M25 }

\end{abstract}

\renewcommand\theequation{{\thesection{}.\arabic{equation}}}
\def\setnumbers{\setcounter{equation}{0}}


\section{Introduction}\label{sec:1}
In search of optimal representatives of given knot classes
Fukuhara \cite{fukuhara_1988} proposed the concept of {\it knot energies}
as functionals defined on the space of knots, providing infinite
energy barriers between different knot types.  This concept was made more precise
later and was investigated by various authors; see e.g. \cite{simon_1996}, \cite{buck-simon_1997},
\cite{sullivan_2002}, and we basically follow here  the definition in the book
of O'Hara \cite[Definition 1.1]{ohara-book_2003}.

Let $\mathcal{C}$ be the class of all closed rectifiable curves $\gamma\subset\R^3$ 
whose one-dimensional Hausdorff measure $\H^1(\gamma)$ is equal to 1.
Moreover we assume that all curves in $\mathcal{C}$ contain
a fixed point, say the origin in $\R^3$, and that all loops in $\mathcal{C}$ are parametrized by arclength defined on the interval $[0,1]$, i.e. $\gamma: [0,1] \to \R^3$ is Lipschitz continuous with $|\gamma'|=1$ and $\gamma(0)=\gamma(1)$. 
Members of $\cC$ will sometimes be referred to as {\it (unit) loops.}

\begin{definition}[\textbf{Knot energy}]\label{def:knot-energy}
Any functional $\E:\cC\to [-\infty,\infty]$ that is finite on all simple smooth
loops $\gamma\in\cC$ with the property that
$\E(\gamma_i)$ tends to $ +\infty$ as $i\to\infty$ on any
sequence of simple loops $\gamma_i\in\cC$ that converge uniformly to a limit
curve with at least one self-intersection, is called {\em self-repulsive} or
{\em charge.} If $\E$ is self-repulsive and bounded from below, it is called
a {\em knot energy.}
\end{definition}

One of the most prominent examples of a knot energy is the {\it M\"obius energy},
introduced by O'Hara \cite{ohara_1991}, and written here with a regularization
slightly
different from O'Hara's \cite[p. 243]{ohara_1991}, and
used, e.g., by Freedman, He, and Wang \cite{FHW}:
$$
\emob(\gamma):=\int_0^1\int_0^1\left\{\frac{1}{|\gamma(s)-\gamma(t)|^2}-
\frac{1}{d_\gamma(s,t)^2}\right\}\,ds\, dt\quad\Fo\gamma\in\cC,
$$
which is non-negative, since the intrinsic distance
$$
d_\gamma(s,t):=\min\{|s-t|,1-|s-t|\}
$$
of the two curve points
$\gamma(s),\gamma(t)$  always dominates the extrinsic Euclidean distance
$|\gamma(s)-\gamma(t)|.$  That this energy is indeed self-repulsive is
proven in \cite[Theorem 1.1]{ohara_1994} and \cite[Lemma 1.2]{FHW}, and one
is lead to the natural question if one can minimize $\emob$ in a fixed knot
class. In view of the direct method in the calculus of variations one would try to
establish uniform bounds on minimizing sequences in appropriate norms to pass
to a converging subsequence with limit hopefully in the same knot class.

However, O'Hara observed in \cite[Theorem 3.1]{ohara_1994}
that knots can pull tight in a convergent sequence of loops with uniformly bounded
M\"obius energy. This {\it pull-tight phenomenon} in a sequence of loops $\gamma_i\in\cC$
is characterized by non-trivially knotted arcs $A_i\subset\gamma_i$ of a fixed knot type
that are contained in open balls
$$
B_i\equiv B_{r_i}(x_i)\subset\R^3,\quad\textnormal{such that $r_i\to 0$ as $i\to \infty;$}
$$
see \cite[Definition 1.3]{ohara-book_2003}. In principle this phenomenon could
result in minimizing sequences for $\emob$ of a fixed knot class, converging to a limit
in a different knot class, and it is actually conjectured by Kusner and Sullivan
\cite{kusner-sullivan-moebius_1998} that this indeed happens.
It was one of the great achievements of Freedman,
He, and Wang in their seminal paper
\cite{FHW} to establish the existence of $\emob$-minimizing knots
but restricted to {\it prime knot classes}, with the help of the
invariance of $\emob$ under M\"obius transformations\footnote{This invariance
proven in \cite[Theorem 2.1]{FHW} gave $\emob$ its name.}.
\begin{definition}\label{def:minimizable-tight}
A knot energy $\E$ is {\em minimizable}\footnote{O'Hara calls this property {\it
minimizer producing}; see \cite[Definition 1.2]{ohara-book_2003}.}
if in each knot class there is at least one
representative in $\cC$ minimizing $\E$ within this knot class. $\E$ is called
{\em tight} if $\E(\gamma_i)$ tends to $+\infty$ on a sequence  $\{\gamma_i\}\subset
\cC$
with a pull-tight phenomenon.
\end{definition}
In that sense, $\emob $ is conjectured to be not minimizable (on composite knots)
since it fails to be tight. Changing the powers in the denominators, and allowing for
powers of the integrand, there arises a whole family of different energies, and it
depends on the ranges of parameters whether or not one finds minimizing knots;
see \cite{ohara_1992b,ohara_1994}, \cite{blatt-reiter}.

The purpose of the present note is to investigate knot-energetic properties of
geometrically defined curvature energies involving
{\it Menger curvature}. The basic building block of these functionals
are elementary geometric quantities like the circumcircle radius
$
R(x,y,z)$ of
three distinct points $x,y,z\in\R^3$, the inverse of which
is sometimes referred to as
Menger\footnote{Coined after Karl Menger who intended to develop a purely metric geometry
\cite{menger_1930}; see also the
monograph \cite{blumenthal-menger_1970}.}
 curvature of $x,y,z$.

 Varying one or several of the points $x,y,z$ along the curve $\gamma$ one obtains
 successive smaller radii, whose values then
 depend on the shape of the curve $\gamma$:
 \begin{equation}\label{radii}
 \varrho[\gamma](x,y):=\inf_{z\in\gamma\atop z\not= x\not= y\not=z}R(x,y,z),
 \quad
 \varrho_G[\gamma](x):=\inf_{y,z\in\gamma\atop z\not= x\not= y\not=z}R(x,y,z),
 \quad
{\triangle}[\gamma]:=\inf_{x,y,z\in\gamma\atop z\not= x\not= y\not=z}R(x,y,z).
 \end{equation}
Repeated integrations over inverse powers of
these radii with respect to the remaining variables lead to the
various {\it Menger curvature energies}
\begin{equation}\label{integral-menger}
\M_p(\gamma):=\int_\gamma\int_\gamma\int_\gamma\frac{d\H^1(x)d\H^1(y)d\H^1(z)}{R(x,y,z)^p},
\end{equation}
\begin{equation}\label{double-integral}
\I_p(\gamma):=\int_\gamma\int_\gamma\frac{d\H^1(x)d\H^1(y)}{\varrho[\gamma](x,y)^p},
\end{equation}
and
\begin{equation}\label{lp-global-curvature}
\U_p(\gamma):=\int_\gamma\frac{d\H^1(x)}{\varrho_G[\gamma](x)^p},
\end{equation}
where the integration is taken with respect to the one-dimensional
Hausdorff-measure $\H^1$. By definition \eqref{radii} of the radii
the energy values on a fixed loop $\gamma\in\cC$ are ordered as
\begin{equation}\label{energy-order}
\M_p(\gamma)\le\I_p(\gamma)\le\U_p(\gamma)\le\frac{1}{\triangle[\gamma]^p}
\Foa p\ge 1
\end{equation}
with the limits
\begin{equation}\label{p-limits}
\lim_{p\to\infty}\M_p^{1/p}(\gamma)=
\lim_{p\to\infty}\I_p^{1/p}(\gamma)=
\lim_{p\to\infty}\U_p^{1/p}(\gamma)=\frac{1}{\triangle[\gamma]},
\end{equation}
and each of the sequences $\{\M_p^{1/p}(\gamma)\}$,
$\{\I_p^{1/p}(\gamma)\}$,
$\{\U_p^{1/p}(\gamma)\}$ is non-decreasing as $p\to\infty$ on
a fixed loop $\gamma\in\cC.$

The idea of looking at minimal radii as in \eqref{radii} goes
back to Gonzalez and Maddocks \cite{GM}, where $\varrho_G[\gamma]$
is introduced as the {\it global radius of curvature} of $\gamma$,
and $\triangle$ stands for the {\it thickness} of the curve,
which is justified by the fact that $\triangle$ equals the
classic normal injectivity radius for smooth curves \cite[Section 3]{GM};
see also \cite[Lemma 3]{GMSvdM}
for the justification in the non-smooth case.
The quotient length/thickness (which equals  $1/\triangle$ on the
class $\cC$ of unit loops)  is called
{\it ropelength} and plays a fundamental role in the search
for {\it ideal knots  and links}; see \cite[Section 5]{GMSvdM},
\cite[Section 2]{CKS02},  and \cite{gonzalez-delallave_2003};
see also \cite{heiko3} and \cite{CFKSW06}.
Some knot-energetic properties of ropelength have been
established (see e.g. \cite[Theorems T4 and 4, Corollary 4.1]{buck-simon_1997}),
and we are going to benefit from that.

Allowing  higher order contact of circles (or spheres) to a
given loop $\gamma\in\cC$ one can define various other radii
as discussed in detail in \cite{gonzalez-maddocks-smutny}.
As a particular example we consider the {\it tangent-point} radius
\begin{equation}\label{tan-point-radius}
r_{\textnormal{tp}}[\gamma](x,y)
\end{equation}
as the radius of the unique circle through $x,y\in\gamma$ that is
tangent to $\gamma$ at the point $x$, which, according to Rademacher's
theorem \cite[Section 3.1, Theorem 1]{EG} on the differentiability of Lipschitz functions,
is defined for almost every $x\in\gamma.$
This leads to the corresponding {\it tangent-point} and
{\it symmetrized tangent-point} energy (as mentioned in
\cite[Section 6]{GM})
\begin{equation}\label{tp-energies}
\E_p(\gamma):=\int_\gamma\int_\gamma\frac{d\H^1(x)d\H^1(y)}{r_{\textnormal{tp}}[\gamma](x,y)^p}, \quad
\E^\textnormal{sym}_p(\gamma):=\int_\gamma\int_\gamma\frac{d\H^1(x)d\H^1(y)}{\Big(r_{\textnormal{tp}}[\gamma](x,y)r_{\textnormal{tp}}[\gamma](y,x
)\Big)^{p/2}}
\end{equation}
to complement our list of Menger curvature energies on $\cC$.

\begin{REMARK} In some of our earlier papers, see e.g. \cite{heiko1}, \cite{svdm-MathZ}, \cite{ssvdm-double,ssvdm-triple},
for technical reasons that are of no relevance here, parametric versions of \eqref{integral-menger}--\eqref{lp-global-curvature} and
\eqref{tp-energies} were considered. In the \emph{supercritical} range of parameters that is considered throughout the present paper,
 finiteness of any of these curvature energies implies that 
 $\gamma([0,1])\in \cC$ is homeomorphic to a circle.  
  Therefore, in virtually all the results below we
assume, without any loss of generality, that $\gamma$ is a simple
closed curve, i.e. $\gamma\colon [0,1)\to \R^3$ is injective and
$\gamma(0)=\gamma(1)$.
\end{REMARK}

Why do we care about these energies if there are already
O'Hara's potential
energies such as $\emob$,
and -- as a kind of hard or steric counterpart --
ropelength? First of all, O'Hara's energies require some sort
of regularization due to the singularities of the integrands
on the diagonal of the domain $[0,1]^2$, whereas the
coalescent limit $x,y,z\to\zeta$ on a sufficiently smooth loop $\gamma$
leads to convergence of $1/R$ to classic curvature $\kappa_\gamma(\zeta)$:
$$
\lim_{x,y,z\to\zeta}R^{-1}(x,y,z)=\kappa_\gamma(\zeta),
$$
so that no regularization is necessary
as pointed out by Banavar et al. in \cite{banavar}\footnote{More
on convergence of the various radius functions in \eqref{radii}
in the setting of non-smooth loops can be found in
\cite{heiko1} and \cite{svdm-MathZ}.}.
Moreover, using the elementary geometric definition of the
respective integrands we have gained detailed insight in
the regularizing effects of Menger curvature energies in
a series of papers \cite{svdm-MathZ,ssvdm-double,ssvdm-triple,svdm-tpcurves,kampschulte_2012}. In particular, the uniform $C^{1,\alpha}$-a-priori
estimates
for supercritical values of the power $p$, that is,
for $p$ above the respective critical value, for which the
corresponding energy is scale-invariant,  turn out to be the
essential
tool, not only for compactness arguments that play a central role
in variational applications, but also in the present knot-theoretic
context; see Section 2. Let us mention that even in the subcritical case
these energies may exhibit regularizing behaviour if one starts on
a lower level of regularity, e.g. with measurable sets \cite{leger},\cite{mattila-lin_2001}, \cite{scholtes_2011,scholtes_2012}. Integral
Menger curvature $\M_2$, for example, plays a fundamental role
in harmonic analysis for the solution of the Painlev\'e problem;
see \cite{melnikov-verdera,mattila-melnikov-verdera,verdera_2001,david-survey,tolsa-survey,dudziak_2010}.
Moreover, in contrast to O'Hara's repulsive potentials, the elementary geometric integrands in \eqref{radii} have lead to
higher-dimensional analogues of discrete curvatures
where one can establish similar
$C^{1,\alpha}$-estimates for a priorily
non-smooth admissible sets of finite energy of arbitrary dimension
and co-dimension \cite{StvdM1,StvdM2,svdm-surfaces,svdm-tp1,skol2,
kolasinski_2012,ksvdm-w2p,
ksvdm-tp2},
which could initiate further analysis of higher dimensional knot space.
The problem of finding a higher-dimensional variant of, e.g., the M\"obius energy that is analytically
accessible to variational methods seems wide open; see \cite{kusner-sullivan-moebius_1998,auckly-sadun_1997,fuller-vemuri_2010}.
Finally, recent work of Blatt
and Kolasinski \cite{blatt-tp,blatt-menger_2011},\cite{blatt-k}   characterizes the energy spaces of Menger-type curvatures in terms
of (fractional) Sobolev spaces, so that one can hope to tackle
evolution problems for integral Menger curvature $\M_p$, for instance,
in order to untangle complicated configurations of the unknot, or to
flow complicated representatives of a given knot class to a
simpler configuration without leaving the knot class; see
recent numerical work of Hermes in \cite{hermes-phd}.

In order to investigate knot-energetic properties of the Menger curvature
energies in \eqref{integral-menger}--\eqref{lp-global-curvature} and \eqref{tp-energies} in more depth we will discuss three more
properties (cf. \cite[Definition 1.4]{ohara-book_2003}).
\begin{definition}\label{def:strong-unknot-detecting-basic}
\begin{enumerate}
\item[\rm (i)]
A knot energy $\E$ on $\cC$ is {\em strong} if there are only finitely many
distinct knot types under each energy level.
\item[\rm (ii)]
A knot energy $\E$  {\em distinguishes the unknot} or is called {\em unknot-detecting}
if the infimum of $\E$ over the trivial knots (the ``unknots'') in $\cC$
is strictly
less than the infimum of $\E$ over the non-trivial knots in $\cC$.
\item[\rm (iii)]
A knot energy $\E$ is called {\em basic} if the round circle is the unique
minimizer of $\E$ in $\cC.$
\end{enumerate}
\end{definition}
Many of the knot-energetic properties we establish here for Menger curvature energies
can be summarized in the following table, where for comparison we have included
the M\"obius energy $\emob$ and also {\it total curvature}
\begin{equation}\label{total_curvature}
TK(\gamma):=\int_\gamma |\kappa_\gamma|\,ds\quad\textnormal{for sufficiently smooth
loops $\gamma\in\cC$,}
\end{equation}
even though this energy as an integral over classic curvature, that is, over a
purely local quantity, does not even detect self-intersections, so that total
curvature fails to be a knot energy altogether.


\begin{center}
\begin{tabular}{c||c|c|c|c|c|c||c|c}
Is the energy: & $\M_{p>3}$ & $\I_{p>2}$ & $\U_{p>1}$ & $\E_{p>2}$ & $\E^\textnormal{sym}_{p>2}$
& $1/\triangle$ & $\emob$ &
$TK$ \\
\hline %
charge             & Yes & Yes & Yes & Yes & Yes  & Yes & Yes & No \\
minimizable        & Yes & Yes & Yes & Yes & Yes  & Yes & No  & No \\
tight              & Yes & Yes & Yes & Yes & Yes  & Yes & No  & No \\
strong             & Yes & Yes & Yes & Yes & Yes  & Yes & Yes & No \\
unknot-detecting   & ?   & Yes & Yes & Yes & Yes & Yes  & Yes & Yes \\
basic              & ?   & ?   & Yes & ?   & ?    & Yes & Yes & No \\
       \end{tabular}
       \end{center}
The respective index of each Menger curvature energy in this table
denotes the admissible
supercritical range of the power $p$, where we have neglected the fact that most
of these energies do penalize  self-intersections even in the scale-invariant case,
that is, curves with double points
have infinite energies $\I_2$, $\U_1$, $\E_2$; see \cite[Proposition 2.1]{ssvdm-double},
\cite[Lemma 1]{svdm-MathZ}, and \cite[Theorem 1.1]{svdm-tpcurves}.

Notice that the affirmative answers in the first five columns settle conjectures
of Sullivan \cite[p. 184]{sullivan_2002} and O'Hara \cite[p. 127]{ohara-book_2003}
at least for the respective supercritical range of $p$ and for one-component links.

The M\"obius energy $\emob$ is strong since it bounds the {\it average crossing number}
$\acn$
that according to \cite[Section 3]{FHW} can be written as
\begin{equation}\label{acn}
\acn(\gamma):=\frac{1}{4\pi}\int_0^1\int_0^1\frac{|(\gamma'(s)\times\gamma'(t))\cdot
(\gamma(t)-\gamma(s))|}{|\gamma(t)-\gamma(s)|^3}\,ds\, dt\quad\Fo\gamma\in\cC,
\end{equation}
where $\times$ denotes the usual cross-product in $\R^3.$
As a consequence of the good bound obtained in
\cite[Theorem 3.2]{FHW} Freedman, He, and Wang can show
that $\emob$ also distinguishes the unknot; see \cite[Corollary 3.4]{FHW}.
In \cite{abrams-etal_2003} it is moreover shown that $\emob$ is basic (as well as many other repulsive potentials), which settles the column for $\emob$ in the table above.
The only ``Yes'' for total curvature is due to the famous Far\`y-Milnor theorem
\cite{fary_1949}, \cite{milnor_1950}, which establishes the sharp lower bound
$4\pi$ for the total curvature of non-trivially knotted loops, whereas the round circle has
total curvature $2\pi$.   Fenchel's theorem ascertains  the nontrivial lower
bound $2\pi$ for $TK(\gamma)$ for any continuously differentiable
loop $\gamma$ with equality if and only if $\gamma$ is a planar simple convex curve,
which, however,  does not suffice to single out the circle as the only minimizer, so
$TK$ is not basic.

To justify the affirmative entries in the first four rows for
the Menger curvature energies we are going to use compactness arguments
based on the respective a priori estimates we obtained in our earlier work.
This is carried out in Section 2. The properties ``unknot-detecting''
and ``basic'' are dealt with individually in Section 3, and there are some additional observations.
The great circle on the boundary of a ball uniquely minimizes $\I_p$ for every $p\ge 2$
among all curves packed into that ball (Theorem \ref{thm:3.2}). This restricted
version of the property ``basic'' is accompanied by a non-trivial lower bound for
$\I_p$ (Proposition \ref{prop:3.3}), and the observation that $\gamma$ must
be a circle if $R$, or $\varrho[\gamma]$, or $\varrho_G[\gamma]$ is constant
along $\gamma$. In addition,
we show that any minimizer of integral
Menger curvature $\M_p$ is unknotted if $p$ is sufficiently large\footnote{Ropelength, $\emob$, and
$\U_p$ for all $p\ge 1$  are basic, which, of course, is a much stronger property.}.
 In Section 4
we prove additional properties relevant for knot-theoretic
considerations. In Theorem \ref{thm:4.1} we show that polygons
inscribed in a loop of finite energy and with vertices
spaced
by some negative power of the energy value are isotopic to the
curve. This  produces a bound on the  stick number (Corollary \ref{cor:stick_number}) and therefore
also an alternative direct proof for these energies to be strong (Corollary \ref{cor:finiteness}); cf.
\cite[Theorem 2, Corollary 4]{litherland-etal_1999}
for related results for ropelength.
Theorem \ref{thm:4.1}  can also be used to prove that the energy level
of two loops $\gamma_1,\gamma_2\in\cC$ determines a bound on
the Hausdorff-distance $d_\mathscr{H}(\gamma_1,\gamma_2)$
below which the two curves are isotopic (Theorem \ref{thm:4.2}).
Both results rely on a type of
excluded volume  and restricted bending constraint that finite energy imposes on the curve,
that we refer to as ``diamond property'' (see Definition
\ref{def:4.0}),  which is much weaker than positive reach \cite{federer_1959}. It
does not mean that there is  a uniform neighbourhood guaranteeing the unique
next-point projection, which would correspond to finite ropelength; see \cite[Lemma 3]{GMSvdM}.  Roughly speaking, it means that any chain
of sufficiently  densely spaced points carries along a ``necklace''
of diamond shaped regions as the only permitted zone for the curve
within a larger tube; see Figure~\ref{2cones}.

\subsection{Open problems}\label{sec:1.1}
The question marks in the table above depict unsolved problems. In particular the question if $\M_p$, $\I_p$, $\E_p$, or $\E_p^\textnormal{sym}$
are basic remains to be investigated.
 Notice, however,
 that Hermes recently proved that the circle is a critical point
 of $\M_p$ \cite{hermes-phd}, which also supports our conjecture that all these energies
 are basic. Numerical experiments suggest
that $\M_p$ should clearly distinguish the unknot, but so far we have not been able to  prove that. Our bounds for the average crossing number
$\acn$ are by far not good enough to  capture that.
Moreover,
our compactness arguments to prove the properties in the first four
rows of the table work in the respective supercritical case, i.e.,
for $p>3$ for $\M_p$, for $p>2$ for $\I_p,\E_p$, and $\E_p^\textnormal{sym}$, and for $p>1$ for $\U_p$. But what happens for the geometrically
interesting scale-invariant cases $\M_3,\I_2,\E_2,\E_2^\textnormal{sym},
\U_1$? \heikodetail{\tt\xx It could be mentioned that Blatt's  and Kolasinski's  characterization of the
energy spaces in terms of fractional Sobolev spaces covers part of the subcritical
range, too, since they assume $C^1$-curves from the beginning on, and this observation
could turn out to be of help; see \cite{blatt-tp,blatt-menger_2011,blatt-k}.}

Further open problems include the regularity theory for minimizing knots of these energies (are they just $C^{1,\alpha}$, as all other curves of finite energy, or $C^{1,1}$, as the minimizers of the ropelength functional, cf. \cite[Theorem~7]{CKS02} and \cite[Theorem~4]{GMSvdM}, or maybe $C^\infty$, like the minimizers\footnote{
See \cite{FHW}, \cite{He}, \cite{blatt-reiter-schikorra_2012} for the regularity theory
for minimizers and critical points of $\emob$, and for less geometric energies that
are related to $\E_p$ see the very recent account \cite{blatt-reiter_2012}.}    of $\E_{\text{M\"{o}b}}$?), and better bounds --- sharp for some knot families, if possible ---  for the average crossing number and stick number in terms of $\M_p$ and other energies, especially in the scale invariant cases mentioned above. Even partial answers would enlarge our knowledge of these curvature energies and their global properties.

\heikodetail{\tt \underline{Further open questions:}
\begin{itemize}
\item
better bounds for acn, stick number etc. for $\M_p$?

\item \yy the growth of infima of energies as the complication of knot class grows (say, in terms of crossing number or acn)?
\end{itemize}
}

\medskip\noindent\textbf{Acknowledgement.} The authors gratefully acknowledge: R. Ricca for organizing the ESF Conference on Knots and Links in Pisa in 2011, where parts of this work have been presented in progress;
M. Giaquinta, for repeatedly hosting our stays at the Centro de Giorgi; and R. Kusner for creating the opportunity to finish this work at the Kavli Institute of Theoretical Physics (KITP) in Santa Barbara. This paper grew out of a larger project on geometric curvature energies financed, along with two workshops in B\c{e}dlewo and Steinfeld, by DFG and the Polish Ministry of Science. At KITP this research was supported
in part by the National Science Foundation under Grant No. NSF PHY11-25915.

\section{Charge, strong, and tight}\label{sec:2}
We denote by $C^0([0,1])$ the space of continuous functions and recall the sup-norm
$$
\|f\|_{C^0([0,1])}:=\sup_{s\in [0,1]}|f(s)|\quad\Fo f\in C^0([0,1]),
$$
and the H\"older seminorm
$$
[f]_{0,\alpha}:=\sup_{s,t\in [0,1]\atop s\not=t}\frac{|f(s)-(t)|}{|s-t|^\alpha},
$$
which together with the sup-norm constitutes the $C^{0,\alpha}$-norm
$$
\|f\|_{C^{0,\alpha}([0,1])}:=\|f\|_{C^0([0,1])}+[f]_{0,\alpha}.
$$
The  higher order spaces $C^k([0,1])$, and $C^{k,\alpha}([0,1])$ consist of those
functions that are $k$ times continuously differentiable on $[0,1]$
such that the sup-norm,
respectively the sup-norm
and the H\"older
seminorm of the $k$-th derivative
are finite.

\begin{theorem}\label{Metathm:A}
Let $\E:\cC\to (-\infty, \infty]$ be bounded from below such that
\begin{enumerate}
\item[\rm (i)] 
There exists $\delta = \delta(E)>0$, such that for all curves  
$\gamma \in \cC$ with $\E(\gamma) < E$
\[
|\gamma(s) - \gamma(t)| > \min\big\{\delta, \frac{d_\gamma(s,t)}{2}\big\}.
\]
\item[\rm (ii)]
$\E$ is sequentially
lower semi-continuous on $\cC\cap C^1([0,1],\R^3)$
with respect to $C^1$-convergence.
\item[\rm (iii)] There exist  constants $C=C(E)$ and $\alpha=\alpha(E)
\in (0,1]$ depending only on the energy level
$E$ such that for all $\gamma\in\cC$
with $\E(\gamma)\le E$ one has $\gamma\in C^{1,\alpha}([0,1],\R^3)$ with
 $\|\gamma\|_{C^{1,\alpha}([0,1],\R^3)}\le C.$
\end{enumerate}
Then $\E$ is charge, minimizable, tight, and strong.
\end{theorem}
As an essential tool for the proof of this theorem let us recall that
isotopy type is stable under $C^1$-convergence. In the $C^2$-category one finds this result, e.g., in Hirsch's book \cite[Chapter 8]{hirsch}, whereas the only
published proofs in $C^1$ we are aware of are in the papers by Reiter \cite{reiter_isop_2005} and
by Blatt in higher dimensions \cite{blatt_isop_2009}.
\heikodetail{\tt We could probably prove all that without Reiter's result, since we have
our Theorem \ref{thm:4.2}, but I guess a remark like that would rather
confuse the reader... In addition, in \cite[Prop. 5.1 \& Prof. 5.2]{svdm-tpcurves} we   have proven
that $\E_p$ is charge and strong without using Theorem \ref{thm:isotopy}}
\begin{theorem}[\textbf{Isotopy}]\label{thm:isotopy}
For any curve $\gamma\in C^1([0,1],\R^3)\cap\cC$ there is $\epsilon_\gamma>0$
such that all closed curves $\beta\in C^1([0,1],\R^3)$ with
$\|\beta-\gamma\|_{C^1([0,1],\R^3)}<\epsilon_\gamma$
are ambient isotopic to $\gamma.$
\end{theorem}
{\sc Proof of Theorem \ref{Metathm:A}:}\,
Assume  that $\E$ is not charge, so that one finds a
sequence of  simple curves $\{\gamma_i\}\subset\cC$
with uniformly bounded energy $\E(\gamma_i)\le E<\infty$,
converging uniformly, that is, in the sup-norm
to  $\gamma$ which is not a simple loop.  
By assumption (i) there exist $\delta >0 $ such that for all $i
\in \N$ we have
$|\gamma_i(s) - \gamma_i(t)| \ge \min\{\delta, d_\gamma(s,t)/2\}$. As $\gamma$ is not embedded,
there exist $s\not= t \in [0,1)$ such that $\gamma(s) = \gamma(t)$ and for sufficiently
large $i$ we have
\[
\min\big\{\delta, \frac{d_\gamma(s,t)}{2} \big\} > |\gamma_i(s) - \gamma(s)| + |\gamma_i(t) - \gamma(t)|
\ge |\gamma_i(s) - \gamma_i(t)|,
\]
a contradiction. So, $\E$ is indeed charge.  

Now we would like to minimize $\E$ on a given knot class $[K]$ within $\cC$. Note
first that by rescaling a smooth and regular representative of $[K]$
to length one and reparametrizing to arclength, we find that there is a representative of $[K]$ in $\cC$.
In particular, there is
a minimal sequence $\{\gamma_i\}\subset\cC$ with $\gamma_i\in [K]$ for all $i\in\N$, such that
$$
\lim_{i\to\infty}\E(\gamma_i)=\inf_{\cC\cap [K]}\E,
$$
and the right-hand side is finite, since by assumption $\E$ is bounded from below.
 Thus the sequence of energy values $\E(\gamma_i)$
is uniformly bounded by some constant $E\ge 0$ for all $i\in\N$ and by assumption (iii)  there exist
 constants $C=C(E)$ and $\alpha=\alpha(E)\in (0,1]$
 depending only on $E$ but not on $i\in\N$, such that
$$
\|\gamma_i\|_{C^{1,\alpha}([0,1],\R^3)}\le C(E).
$$
Thus, this sequence is equicontinuous, and by the theorem of Arzela-Ascoli
we can extract a subsequence $\{\gamma_{i_k}\}\subset\{\gamma_i\}$ such that
$\gamma_{i_k}$ converges to $\gamma$ in $C^1$, so that in particular $|\gamma'|
\equiv 1$.
  Assumption (i) implies that all $\gamma_{i_k}$ in
the sequence are simple and,  as  $\E$ is charge, the limit curve is injective, hence
$$
\H^1(\gamma)=\int_0^1|\gamma'(s)|\,ds=
\lim_{k\to\infty}\int_0^1|\gamma_{i_k}'(s)|\,ds=\lim_{k\to\infty}\mathscr{H}^1(
\gamma_{i_k})=1
$$
because of the continuity of the length\footnote{Length is only lower semicontinuous with respect
to uniform convergence, so that a priori $\gamma$ could have had length smaller
than one.}  functional
$$
\textnormal{length}(\gamma):=\int_0^1|\gamma'(s)|\,ds
$$
with respect to $C^1$-convergence. Therefore the limit curve $\gamma$ is
in $\cC\cap C^1([0,1],\R^3).$
 We can use assumption (ii) to conclude
$$
\E(\gamma)\le\liminf_{k\to\infty}\E(\gamma_{i_k})=\inf_{\cC\cap [K]}\E.
$$

According to Theorem \ref{thm:isotopy} we find that
$$
[K]=[\gamma_{i_k}]=[\gamma]\quad\Foa k\gg 1,
$$
so that $\gamma\in\cC\cap [K]$ and therefore
$$
\inf_{\cC\cap [K]}\E\le \E(\gamma)\le \inf_{\cC\cap [K]}\E,
$$
i.e. equality here, which establishes $\gamma\in\cC$ as the
(in general not unique) minimizer.

As to proving that $\E$ is tight  we assume that there is a sequence with
the pull-tight phenomenon with uniformly bounded energy. As above we find
a $C^1$-convergent subsequence $\gamma_{i_k}$ with a $C^1$-limit curve $\gamma\in\cC$
that necessarily has the same knot type as $\gamma_{i_k}$ for all sufficiently large
$k$ according to Theorem \ref{thm:isotopy}. But this contradicts the fact that
a subknot is pulled tight which would change the knot-type in the limit.
\heikodetail{\tt a bit of hand-waving here, which arises from the definition of pull-tight
phenomenon which is not so precise...\yy OK for me, P.} Consequently, $\E$ is tight.

Assume finally that there are infinitely many knot-types $[K_i]$
with representatives
$\gamma_i\in\cC$  with uniformly bounded energy $\E(\gamma_i)\le E$ for all $i\in\N.$
Again, we extract a subsequence $\gamma_{i_k}\to \gamma$ in the $C^1$-topology.
Hence $\E(\gamma)\le E$ by assumption (ii), so that $\gamma$ is embedded
by assumption (i). But then by means of Theorem \ref{thm:isotopy} we reach
a contradiction to infinitely many knot-types by
$
[\gamma]=[\gamma_{i_k}]$ for all sufficiently large $k\in\N.$
Consequently, $\E$ is also strong, which concludes the proof of the
theorem.
\qed

\begin{corollary}\label{cor:2.2}
The energies $\M_p$ for $p>3$, $\E_p$, $\esym_p$, and $\I_p$ for $p>2$, and
$\U_p$ for $p>1$ are charge, minimizable, tight, and strong.
\end{corollary}
\proof
We need to check the validity of the assumptions in Theorem \ref{Metathm:A}
for each of the energies under consideration.

For $\I_p$ assumption (i) follows from  \cite[Prop. 3.3]{ssvdm-double}, (ii) from \cite[Lemma 3.5]{ssvdm-double}, and (iii) from Corollary 3.2 of \cite{ssvdm-double}.

For $\U_p$, assumption (ii) of Theorem~\ref{Metathm:A} follows from
\cite[Thm.~3]{svdm-MathZ}, (iii) of Theorem~\ref{Metathm:A}
is provided by \cite[Thm.~1~(iv)]{svdm-MathZ}, 
and to verify (i) one can use the results from 
Chapter 4 of the present paper, namely Lemma \ref{34-bilip} and Proposition \ref{prop diam}.\footnote{Another idea to check that $\U_p$ is charge, minimizable, tight and strong is to notice that the
statement of Theorem~\ref{Metathm:A} holds also 
if instead of (i) we assume that every curve with finite 
energy is injective  and that the energy is charge 
-- $\U_p$ satisfies these requirements by Lemma~1 and 
Theorem~3 from \cite{svdm-MathZ}.  }


For integral Menger curvature $\M_p$ the situation  is a little more subtle.
By Theorem~1.4 in \cite{ssvdm-triple}, if $\M_p(\gamma)<\infty$, then $\gamma([0,1])$ is, topologically, a segment or a circle.  By the very definition of $\cC$ for each $\gamma\in \cC$ the first possibility can easily be excluded and we deal in fact with an arclength parametrized, simple closed curve $\gamma$. Thus, one can apply Theorems 1.2 and 4.3 in \cite{ssvdm-triple} to obtain the a priori estimate needed in (iii) along with (i), whereas (ii) is dealt with in Remark 4.5 of that paper.

To justify (i) for $\E_p$ we need to combine
Theorem 1.1 in \cite{svdm-tpcurves} with the aforementioned
result of \cite{ssvdm-triple}, which requires only one simple point of the locally homeomorphic curve to deduce injectivity of the arc\-length parametrization. (One simply has to copy the arguments in \cite[Section 3.1]{ssvdm-triple} to extend the proof of Theorem
3.7 from that paper to cover the case of the tangent-point energy $\E_p$.) Uniform $C^{1,\alpha}$ bounds are given in Proposition~4.1 of \cite{svdm-tpcurves}. With that information, the verification of lower semicontinuity of $\E_p$ on $\cC\cap C^1$ is a simple exercise, requiring an application of Fatou's lemma. Indeed,  since 
$$\frac{1}{r_{\text{tp}}(\gamma(t),\gamma(s))}=\frac{2\dist(\gamma(s),\ell(t))}{|\gamma(t)-\gamma(s)|^2},$$
where $\ell(t)=\{\gamma(t)+\tau \gamma'(t)\, \colon \, \tau\in\R\}$ is the tangent line to $\gamma$ at $\gamma(t)$, for curves $\gamma_j\to \gamma$ in $C^1$ we obviously have $1/r_{\text{tp}} (\gamma_j(t),\gamma_j(s))\to 1/r_{\text{tp}} (\gamma(t),\gamma(s))$ whenever the limit is nonzero, and the result follows.

As to $\esym_p$ we refer to Lemma 3.1.3  in \cite{kampschulte_2012}, and in particular to the final
estimate of its proof in \cite[p. 23]{kampschulte_2012}
to verify assumption (iii). Assumptions (i) and (ii) can be proven as we indicated above for $\E_p$.
\qed



\heikodetail{

{\tt\yy The old material 2.1 is in the file, hidden after \verb#\end{document}#}

{\tt Acknowledged!}

}

\section{Basic, detecting unknots}\label{sec:3}

\setnumbers

\begin{proposition} The energy $\U_p$ is basic and unknot-detecting for all $p\ge 1$.
\end{proposition}

\begin{proof}
We may assume $\U_p(\gamma)<\infty$, so that by \cite[Theorem 1]{svdm-MathZ} $\gamma$ is in the Sobolev space $W^{2,1}$ of twice weakly
differentiable functions with weak second derivatives in $L^1$. In particular, classic local curvature $\kappa_\gamma$ of $\gamma$
exists almost everywhere.    Since $\varrho_{G}$ does not exceed the
local radius of curvature wherever the latter exists \cite[Lemma 7]{heiko1}, we can estimate
\[
4\pi \le \int_\gamma |\kappa_\gamma|\, ds \le \left(\int_\gamma
|\kappa_\gamma|^p\, ds\right)^{1/p} \le \left(\int_\gamma
\big({\varrho_G}[\gamma]\big)^{-p}\, ds\right)^{1/p} = \U_p^{1/p}(\gamma)
\]
 for any  non-trivially knotted curve $\gamma$ of
length 1
by the  Far\'{y}--Milnor theorem and H\"{o}lder's inequality,
whereas
\[
\U_p^{1/p}(\text{circle}) =2\pi < 4\pi\, .
\]
 Lemma 7 in \cite{svdm-MathZ} states that the circle uniquely minimizes $\U_p$.   \qed
\end{proof}

We do not know if the energies $\M_p$, $\I_p,$ $\E_p$, and $\esym_p $  are basic or not. But for $\I_p$ we can
at least prove a restricted version of that property, which may be interpreted as a relation between energy and compaction:
When stuffing a unit loop into a closed
ball the most energy efficient way (with respect to $\I_p$) is to form a great circle. Buck and
Simon have established a non-trivial lower bound for their {\it normal energy}  for curves packed into a ball  in  \cite[Theorem 1]{buck-simon_1997},
however, without presenting an explicit  minimizer.
It turns out that this normal energy is  proportional to the tangent-point energy $\E_2$, and one
might hope to use their  bound for $\I_2$ by the simple ordering $\I_2\ge \E_2$ (cf. \eqref{ineq2} in the proof of
Corollary \ref{cor:3.5} below). But we obtain a better bound  for $\I_p$ using a powerful sweeping argument which requires
the infimum in the definition of the particular radius $\varrho[\gamma]$ in the integrand. Moreover, this technique of proof permits our uniqueness argument.
\begin{theorem}[\textbf{Optimal packing in ball}]\label{thm:3.2}
Among all loops in $\cC$ that are contained in a fixed closed ball $\overline{B_{\frac{1}{2\pi}}}$,
circles of length $1$, i.e.,  great circles on $\partial B_{\frac{1}{2\pi}}$, ``uniquely'' minimize $\I_p$
for all $p\ge 2.$
\end{theorem}
We start with a technical lemma that also contains the aforementioned sweeping argument.
\begin{lemma}[\textbf{Sweeping}]\label{lem:sweeping}
Let $\gamma\subset\cC$, and assume that there are two distinct points $x,y\in\gamma$  with
\begin{equation}\label{sweeping-cond}
\rho:=\varrho[\gamma](x,y)>\frac{|x-y|}{2}.
\end{equation}
\heikodetail{\tt since so far we were using $\rho$ and $\varrho$ in an arbitrary way, I tried
to adapt to $\varrho$ whenever we speak of one of our special radii, for intermediate
definitions or abbreviations I stuck with $\rho$ . If you spot more of those previous
inconsistencies, go ahead...}
Then no point of $\gamma$ is contained in the ``sweep-out region''
\begin{equation}\label{sweep-out_region}
S(x,y):=\bigcup_{x,y\in\partial B_\rho} B_\rho\setminus\overline{\bigcap_{x,y\in\partial B_\rho} B_\rho},
\end{equation}
which is the union of all balls of radius $\rho $ containing $x$ and $y$ in their boundary $\partial B_\rho$ minus
the closure of their intersection.

If, moreover, $|x-y|<\diam\gamma$, or if the weaker assumption  $|x-y|\le |\xi-\eta|$ for at least one pair $(\xi,\eta)\in\gamma\times\gamma\setminus \{(x,y)\}$
holds,
 then $\gamma$ is not completely contained in the lens-shaped
region
\begin{equation}\label{lense}
\ell(x,y):=\overline{\bigcap_{x,y\in\partial B_\rho} B_\rho},
\end{equation}
and we have the estimate
\begin{equation}\label{sharp-estimate}
1\ge |x-y|+\Big(2\pi-2\arcsin \frac{|x-y|}{2\rho}\Big)\rho,
\end{equation}
in particular,
\begin{equation}\label{rough-estimate}
\rho\le\frac{1}{\pi+2}.
\end{equation}
\end{lemma}
\proof
The first statement follows from elementary geometry.
\begin{figure}[!h]
\centering
\includegraphics[height=5cm]{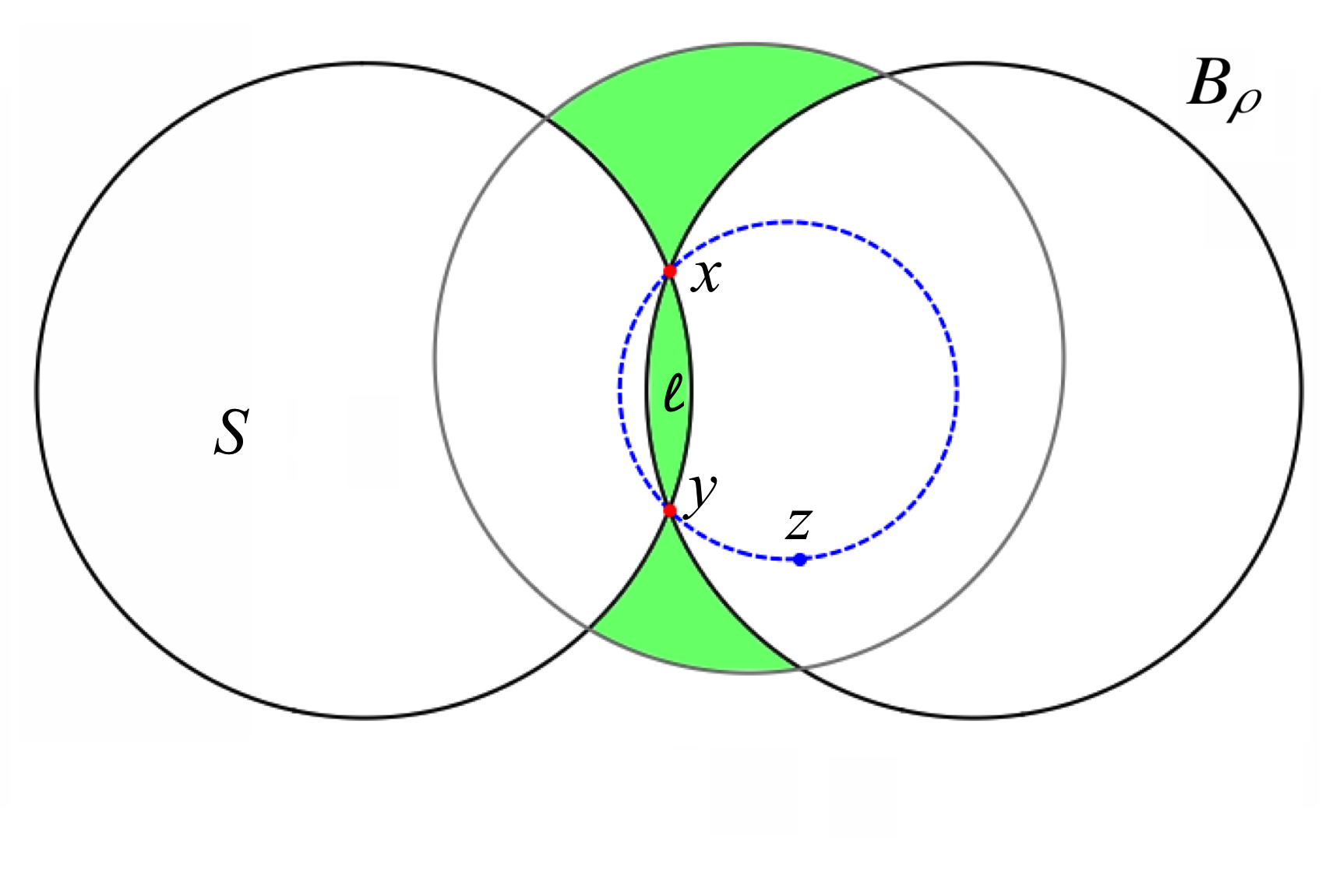}\hfill
\includegraphics[height=5.9cm]{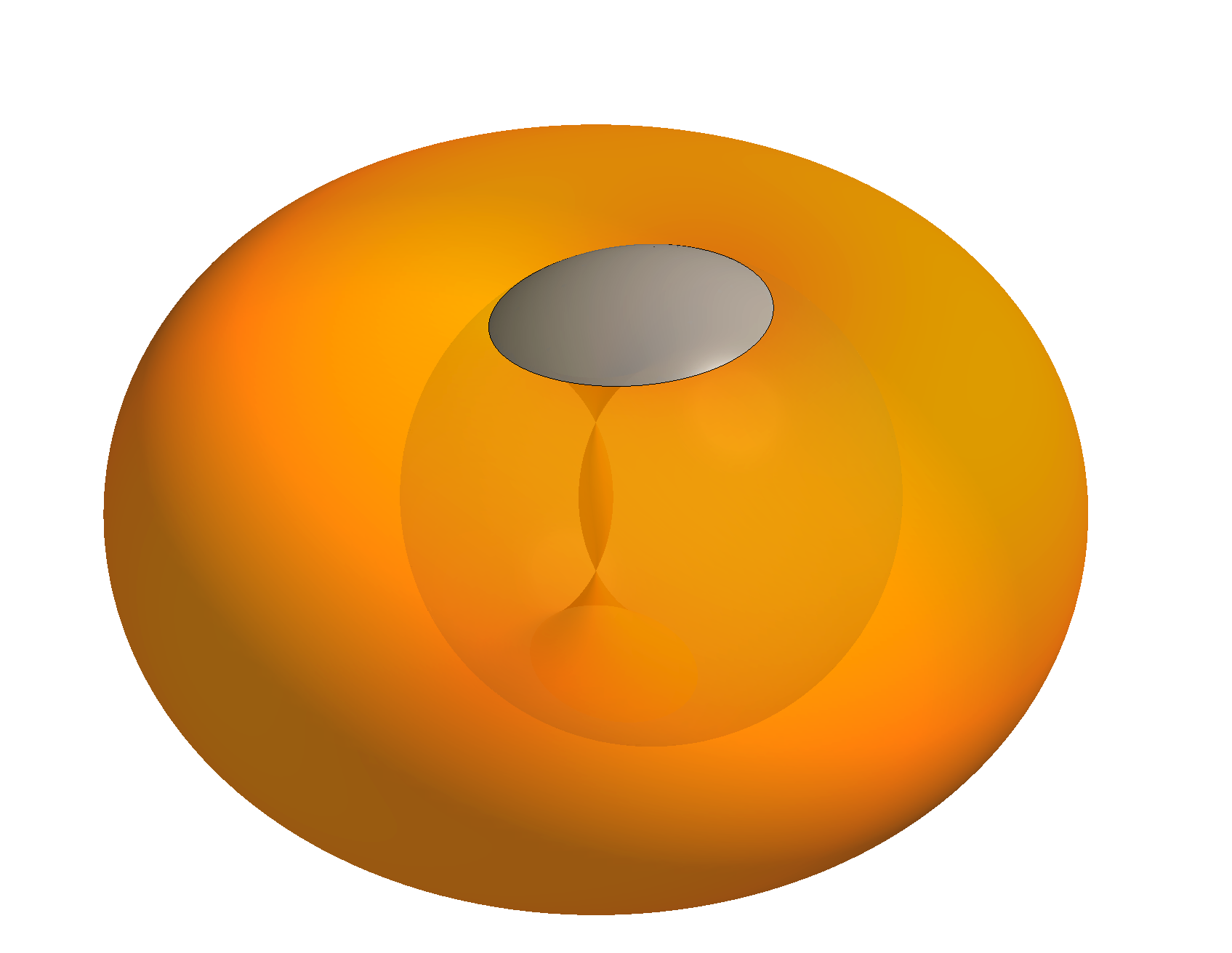} \hfill
$\phantom{l}$
\caption{\label{sweepregion}
Left: If $x,y\in \gamma$ and $\varrho[\gamma](x,y)>|x-y|/2$, then $\gamma\cap S=\emptyset$ for a fairly large region $S=S(x,y)$. In the situation
of Theorem \ref{thm:3.2}  $\gamma$ is confined to the shaded zones.  The solid circle in the middle depicts the boundary of the ball $B_{1/2\pi}$ containing $\gamma$.
Right: A three-dimensional view of the sweep-out region $S$, whose boundary coincides with a self-intersecting torus of rotation, and of the ball $B_{1/2\pi}$, of which a substantial portion is immersed in that torus.}
\end{figure}
Indeed, if there were $z\in S(x,y)\cap\gamma$, then by elementary geometry carried out in the plane spanned by $x,y,$ and $z$,
we would find
$$
R(x,y,z)<\rho=\varrho[\gamma](x,y),
$$
contradicting the very definition of $\varrho[\gamma](x,y)$; see the the dashed circle
with radius $R(x,y,z)$ in the left image of Figure \ref{sweepregion}.

As to the second statement we use the weaker assumption and suppose to the contrary that $\gamma\subset
\ell(x,y)$ which implies a direct contradiction via
$$
\diam\gamma=|x-y| > |\xi-\eta|\quad\Foa (\xi,\eta)\in\gamma\times\gamma\setminus\{(x,y)\}.
$$
Next, observe that since $\gamma$ is a simple closed curve connecting $x$ and $y$, its unit  length is bounded from below by
the shortest possible  simple loop connecting $x$ and $y$ without staying in $\ell(x,y)$ and without entering $S(x,y)$. Such a  loop is the straight segment from $x$ to $y$ together with the circular great arc on the boundary of one of the balls $B_\rho$; hence \eqref{sharp-estimate}
follows. The rough estimate \eqref{rough-estimate} stems from comparing to
 the worst case scenario, when $x$ and $y$ are antipodal on a ball $B_\rho.$
\qed

{\sc Proof of Theorem \ref{thm:3.2}.}\,
We will simply say ``circle'' when we refer to a circle of unit length, i.e., a round circle in $\cC$.
It suffices to prove the statement for $p=2$, since
\begin{equation}\label{circle_strict}
\I_2(\text{circle})<\I_2(\gamma)
\end{equation}
 for any $\gamma\in\cC$ different
from the circle implies by H\"older's inequality
$$
\I_p^{1/p}(\text{circle})=2\pi=\I_2^{1/2}(\circle)<\I_2^{1/2}(\gamma)\le\I_p^{1/p}(\gamma)
$$
for any $p>2.$

\def\Mplus{M^+(\gamma)}

To show \eqref{circle_strict} we consider first the set $\Mplus$ of pairs of points $x,y\in\gamma$ such that
\begin{equation}\label{bad_points}
\varrho[\gamma](x,y)> \frac{1}{2\pi}=\frac{\diam B_{\frac{1}{2\pi}}}{2}.
\end{equation}
We shall prove that $\Mplus$ contains at most one such pair. If $\Mplus$ is empty, there is nothing to prove. Assume the contrary.
Observe that $\varrho[\gamma](x,y)>|x-y|/2$ for each $(x,y)\in \Mplus$
since $\gamma\subset\overline{B_{\frac{1}{2\pi}}}$, and we can apply the first part of
Lemma \ref{lem:sweeping} to deduce that $\gamma$ has no point in common with
the sweep-out region  $S(x,y)$ defined in \eqref{sweep-out_region}. Next, there can be at most one pair $(x,y)\in \Mplus$ such that
$|x-y|=\diam\gamma$, since $\gamma\cap S(x,y)=\emptyset$ and so $\gamma\subset\ell(x,y),$ which implies
that $|\xi-\eta|<\diam\gamma$ for all pairs $(\xi,\eta)\in\gamma\times\gamma$ different from $(x,y)$.

Now, fix $(x,y)\in \Mplus$ (with $|x-y|=\diam\gamma$, if such a pair exists in $\Mplus$, and arbitrary otherwise).
We claim that there cannot be any other point contained in $\Mplus$.  Indeed,  if there were $(u,w)\in \Mplus\setminus
\{(x,y)\}$, then $|u-w|<\diam\gamma$ and we could apply the second statement of Lemma \ref{lem:sweeping}
to the pair $(u,w)$ replacing $(x,y)$
to conclude that $\gamma\not\subset\ell(u,w)$.  But then the \emph{simple\/} curve
$\gamma$ could not connect the points $u$ and $w$ and remain closed, since the complement
 $B_{\frac{1}{2\pi}}\setminus \big(S(u,w)\cup\ell(u,w)\big)$ is disconnected for each $(u,w)\in \Mplus$;
see Figure \ref{sweepregion}, again with $(u,w)$ replacing $(x,y)$.

Thus, $\Mplus$ contains at most one pair of points. In other words,
\begin{equation}\label{good_points}
\varrho[\gamma](\xi,\eta)\le\frac{1}{2\pi} \quad\textnormal{for all $(\xi,\eta)\in\gamma\times\gamma\setminus (x,y)$,}
\end{equation}
which\footnote{If there is no pair $(x,y)\in\gamma\times\gamma$ satisfying \eqref{bad_points} we find
\eqref{good_points} even  for all $\xi,\eta\in\gamma.$} immediately implies the energy inequality
\begin{equation}
\label{energy_inequality}
\I_2^{1/2}(\gamma)\ge 2\pi=\I_2^{1/2}(\circle).
\end{equation}

To prove uniqueness of the minimizer we assume equality in \eqref{energy_inequality}, which implies by means of
\eqref{good_points} that equality holds in \eqref{good_points} for almost all pairs $(\xi,\eta)\in\gamma\times\gamma.$ Now we claim that the set
$$
M_\textnormal{int}(\gamma):=\{(\xi,\eta)\in\gamma\times\gamma:\varrho[\gamma](\xi,\eta)=\frac{1}{2\pi} \,\,\textnormal{
and $\xi$ or $\eta$ lie in the open ball $B_{\frac{1}{2\pi}}$}\}
$$
contains at most one element. Indeed, for all pairs $(\xi,\eta)\in M_\textnormal{int}(\gamma)$ one
has $\varrho[\gamma](\xi,\eta)>|\xi-\eta|/2$. Assuming that $M_\textnormal{int}(\gamma)$ has at least two elements
we can select $(x,y),(\xi,\eta)$ in that set such that $|x-y|\le|\xi-\eta|$ and apply
 Lemma \ref{lem:sweeping} to obtain $\gamma\cap S(x,y)=\emptyset$ as well as  $\gamma\not\subset\ell(x,y)$ which again implies a contradiction
since $\gamma$ cannot connect $x$ and $y$ within $\overline{B_{\frac{1}{2\pi}}}$.

 So we have shown that almost all $(\xi,\eta)\in\gamma\times\gamma$ satisfy equality in \eqref{good_points}
and
 $\xi,\eta\in\partial B_{\frac{1}{2\pi}}(0)$. If there was any point $z\in\gamma\cap B_{\frac{1}{2\pi}}(0)$ then
a whole subarc $\alpha\subset\gamma$ of positive length would lie in the open ball. Thus,
$\alpha\times\alpha\subset M_\textnormal{int}(\gamma)$, contradicting the statement we just made. Hence $ \gamma$ is completely contained in the boundary $\partial B_{\frac{1}{2\pi}}(0),$
and thus any three points $x,y,z\in\gamma$ must span an equatorial plane, otherwise $R(x,y,z)<1/(2\pi)$. But then
there can be at most one such equatorial plane, which implies that $\gamma$ equals the great circle in that plane.
\qed

The sweeping argument demonstrated in the proof of Theorem \ref{thm:3.2} can also be used to derive the following non-trivial lower bound, which states that one needs at least $\I_2$-energy level $16$ to close up a curve.\footnote{That one needs at least $\I_2$-energy $8$ to close a curve
 can already be shown  directly using the fact that any closed curve of length one is contained in a closed  ball
of radius $1/4$; see Nitsche's short proof in \cite{nitsche_1971}. Now the aforementioned  packing result in  \cite[Theorem 1]{buck-simon_1997}
turns out useful, since $\I_2\ge\E_2$ and the latter is four times their normal energy.}
\begin{proposition}[\textbf{Lower bound for $\mathbf{\I_p}$}]\label{prop:3.3}
For any loop $\gamma\in\cC$ and $p\ge 2$ one has the energy estimate
\begin{equation}\label{I-lower-bound}
\I_p^{1/p}(\gamma)\ge\I_2^{1/2}(\gamma)\ge \min\Big\{2+\pi,\frac{2}{\diam\gamma}\Big\}\ge 4.
\end{equation}
\end{proposition}
\proof
The first inequality is just H\"older's inequality, the last can be seen directly, since the diameter of $\gamma$ is bounded by half of its length.
The  second inequality in \eqref{I-lower-bound}, however, requires a  proof.

First we claim that there is a set $T\subset\gamma\times\gamma$ of positive measure such that for each pair of points $(x,y)\in T$ one has
\begin{equation}\label{tschebycheff}
\frac{1}{\varrho[\gamma](x,y)^2}\le\I_2(\gamma),
\end{equation}
since otherwise we could integrate the reverse inequality to get the contradictory statement
$$
\I_2(\gamma)=\int_\gamma\int_\gamma\frac{d\H^1(\xi)d\H^1(\eta)}{\varrho[\gamma](\xi,\eta)^2}>\I_2(\gamma).
$$
If one  pair $(x,y)\in T$  satisfies $\varrho[\gamma](x,y)=|x-y|/2$ which is bounded from above by $\diam\gamma/2$, then
we obtain
from \eqref{tschebycheff}
$$
\frac{1}{\I_2^{1/2}(\gamma)}\le\varrho[\gamma](x,y)\le\frac{\diam\gamma}{2},
$$
which gives the second alternative  of the minimum in \eqref{I-lower-bound}.

In the other case, $\varrho[\gamma](x,y)>|x-y|/2$ for all $(x,y)\in T$, and  we can apply Lemma \ref{lem:sweeping}
again since we can pick two pairs $(x,y),(\xi,\eta)\in T$ with $|x-y|\le|\xi-\eta|$. This results in
$\gamma\cap S=\emptyset$ and $\gamma\not\subset \ell(x,y)$, where $S$ and $\ell$ are defined in \eqref{sweep-out_region}
and \eqref{lense}, respectively. Then we insert the rough estimate \eqref{rough-estimate} into \eqref{tschebycheff} to obtain the remaining alternative
in the desired estimate \eqref{I-lower-bound}.
\qed

%
%

As another immediate consequence of the sweeping technique we observe
that constant $R$, $\varrho[\gamma]$, or $\varrho_G[\gamma]$ allows only for
the circle. Recall that constant classic local curvature does not imply anything
like that; see, e.g. the construction of arbitrary $C^2$-knots of constant
curvature in \cite{mcatee_2004}.
\heikodetail{

\bigskip

 {\tt check if published in a journal somewhere,
she was a student of Colin Adams, there is another paper by Koch and Engelhardt,
that we could add...?}

}
\begin{corollary}[\textbf{Rigidity}]\label{cor:rigidity}
If there is $R_0\in (0,\infty)$, such that a curve $\gamma\in\cC$ satisfies
either $R(x,y,z)\equiv R_0$, or $\varrho[\gamma](x,y)\equiv R_0$, or $\varrho_G[\gamma](x)\equiv R_0$ for all $x,y,z\in\gamma$, then $R_0=1/(2\pi)$ and $\gamma$ is the
round circle of radius $1/(2\pi)$.
\end{corollary}
\proof
By definition \eqref{radii}    of the respective radii it suffices to prove the
statement under the assumption that $\varrho_G[\gamma](x)=R_0$
for all $x\in\gamma.$ Notice first that this implies $1/\triangle[\gamma]=1/R_0<\infty$
such that by \cite[Lemmata 1 \& 2]{GMSvdM}  $\gamma$ is simple
and of class $C^{1,1}([0,1],\R^3)$. Moreover, by the  elementary geometric
expression for $R$ one finds $\diam\gamma\le 2R_0$.

We claim that, in fact, $\diam\gamma=2R_0$, since if not, we could find
points $x,y\in\gamma$ with
$$
\diam\gamma=|x-y|<2R_0=2\varrho_G[\gamma](x),
$$
and we deduce from the first part of Lemma \ref{lem:sweeping} that $S(x,y)\cap\gamma=\emptyset$, where $S(x,y)$ is
the sweep-out region defined in \eqref{sweep-out_region} for
$\rho:=\varrho_G[\gamma](x)  =R_0\le \varrho[\gamma](x,y).$
Since $x$ and $y$ realize the diameter of $\gamma$ we conclude
that $\gamma$ is completely contained in the lens-shaped
region $\ell(x,y)$ defined in \eqref{lense} for $\rho=\varrho_G[\gamma](x),$
which immediately gives a contradiction, since $\gamma$ is of class
$C^1$ and can therefore have no corner points at  $x$ and $y$.
This proves $\diam\gamma=2R_0$, so that $\gamma$ is contained
in the closure of the ball $B^*:=B_{R_0}(\frac{x+y}{2})$,
since any point on $\gamma$ but
outside the closed slab of width $|x-y|$ and orthogonal to the segment $x-y$ would
lead to a larger diameter, and any point $\zeta\in\gamma\setminus\overline{B^*}$
inside the slab
would lead to the contradiction $R(x,y,\zeta)<R_0=\varrho_G[\gamma](x).$
But with $\gamma\subset\overline{B^*}$ we can apply our best packing
result, Theorem \ref{thm:3.2}, to conclude that $\gamma$ must coincide with a great circle on the boundary $\partial B^*$ because of the identity
$$
\frac{1}{R_0}=\I_2^{1/2}(\gamma)=\I_2^{1/2}(\textnormal{great circle
on $\partial B^*$})=\frac{1}{R_0}.
$$
Since $\gamma\in\cC$ has
length one, we compute $R_0=1/(2\pi).$
\qed

\medskip

The energies $\I_p$, $\E_p$ and $\E_p^{\text{sym}}$ for $p\ge 2$ are also
unknot-detecting. This follows  via simple applications of H\"{o}lder and Young inequalities from a key ingredient which is an inequality, due to Simon and Buck, cf. \cite[Theorem 3]{buck-simon_1997},  between  $\E_2^{\text{sym}}$ and the average crossing number, defined in \eqref{acn}.

Here is the result, for which we present here a short proof for the sake of completeness.

\begin{theorem}[\textbf{Buck, Simon}]  Let $\gamma\in {\mathcal C}$ be a simple curve of class $C^1$. Then
\begin{equation}
    \E_2^{\text{\rm sym}}(\gamma)\ge 16\pi\, \text{\rm acn}\, (\gamma).
  \label{BS inequality}
\end{equation}
\label{thm:BS}
\end{theorem}

\medskip

\begin{proof} The theorem follows from a pointwise inequality between the integrands. To see that, let $s\not=t\in [0,1]$ and $r(s,t)=\frac{\gamma(s)-\gamma(t)}{|\gamma(s)-\gamma(t)|}$.  Set
 \[
\alpha=\alpha(s,t)=\ang\big(\gamma'(s),r(s,t)\big)\, , \qquad  \beta=\beta(s,t)=\ang\big(\gamma'(t),r(s,t)\big)\, ,
\]
and rewrite \eqref{acn} as
\begin{equation}
    \label{acn 2}
    \text{acn}\, (\gamma)=\frac{1}{4\pi}
    \iint_{[0,1]\times [0,1]}
    \frac{\big|\det (\gamma'(s),\gamma'(t),r(s,t))\big|}{|\gamma(s)-\gamma(t)|^2}\; ds\, dt.
\end{equation}
We have
\[
\big|\det (\gamma'(s),\gamma'(t),r(s,t))\big|=|\gamma'(s)\times r(s,t)| \cdot \dist\big(\gamma'(t), \text{span}(\gamma'(s), r(s,t))\big) = \sin \alpha \cdot \sin \varphi,
\]
where $\varphi=\varphi(s,t)$ denotes the angle between $\gamma'(t)$ and  $\text{span}(\gamma'(s), r(s,t))$. Denoting the orthogonal projection of $\R^3$ onto $P:= \text{span}(\gamma'(s), r(s,t))$ by $\pi$, one clearly obtains $$\sin \varphi = |\gamma'(t)-\pi(\gamma'(t))|=\dist\big(\gamma'(t),P\big)\le \big|\gamma'(t)-\langle \gamma'(t),r(s,t)\rangle\, r(s,t)\big|=\sin\beta.$$
Thus,
\begin{equation}
    \label{integrand est}  \frac{\sin\alpha(s,t)\; \sin\beta(s,t) }{|\gamma(s)-\gamma(t)|^2}  \ge   \frac{\big|\det (\gamma'(s),\gamma'(t),r(s,t))\big|}{|\gamma(s)-\gamma(t)|^2} \, .
\end{equation}
The left--hand side above is directly related to the tangent--point radius, as
 a simple geometric argument shows that
\[
\frac{1}{r_{\textnormal{tp}}[\gamma](\gamma(s),\gamma(t)) }=\frac{2\dist(\gamma(t), \gamma(s)+\text{span}\, \gamma'(s))}{|\gamma(s)-\gamma(t)|^2}=\frac{2\sin\alpha(s,t) }{|\gamma(s)-\gamma(t)|}\, .
\]
Hence, \eqref{integrand est} translates to
\[
\frac{1}{r_{\textnormal{tp}}[\gamma](\gamma(s),\gamma(t))} \cdot \frac{1}{r_{\textnormal{tp}}[\gamma](\gamma(t),\gamma(s))} \ge 4\frac{\big|\det (\gamma'(s),\gamma'(t),r(s,t))\big|}{|\gamma(s)-\gamma(t)|^2}   \, .
\]
Integrating, we obtain \eqref{BS inequality}. \qed
\end{proof}

\begin{corollary}[\textbf{Unknot-detecting}]\label{cor:3.5}
 The energies $\I_p$, $\E_p$ and $\E_p^{\text{sym}}$ are unknot-detecting for each $p\ge 2$.
\end{corollary}

\begin{proof}  By H\"older's inequality, for curves of unit length we have
\begin{equation}
    \label{ineq1}
    \F_p(\gamma)^{2/p}\ge \F_2(\gamma), \qquad p\ge 2,
\end{equation}
for each energy $\F_p\in \{\E_p,\E_p^{\text{sym}},\I_p\}$. Besides,
\begin{equation}
    \label{ineq2}
    \I_2\ge \E_2\ge \E_2^{\text{sym}}.
\end{equation}
To verify the second inequality in \eqref{ineq2}, just note
\[
\frac{1}{2}\biggl(\frac{1}{r_{\textnormal{tp}}[\gamma](\gamma(s),\gamma(t))^2} + \frac{1}{r_{\textnormal{tp}}[\gamma](\gamma(t),\gamma(s))^2} \biggr) \ge
\frac{1}{r_{\textnormal{tp}}[\gamma](\gamma(s),\gamma(t))} \cdot \frac{1}{r_{\textnormal{tp}}[\gamma](\gamma(t),\gamma(s))}
\]
and integrate both sides with respect to $(s,t)\in [0,1]^2$.

 In order to check that $\I_2\ge \E_2$
we use the explicit formula for the tangent-point radius from elementary geometry
$$
r_\textnormal{tp}[\gamma](x,y)=\frac{|x-y|}{2\sin\ang(y-x,t_x)},
$$
where we assumed that the unit tangent $t_x$  of $\gamma$ at the point $x$ exists, to express the
denominator in terms of the cross-product of $t_x$ and the unit vector $(x-y)/|x-y|$ to obtain
\begin{equation}\label{rad_comparison}
r_\textnormal{tp}[\gamma](x,y)=\lim_{\gamma\ni z\to x}\frac{|x-y|}{2\Big|\frac{y-z}{|y-z|}\times\frac{x-z}{|x-z|}\Big|}=\lim_{\gamma\ni z\to x}R(x,y,z)
 \ge  \inf_{z\in \gamma} R(x,y,z)=\varrho[\gamma](x,y).
\end{equation}

\noindent
Thus, combining \eqref{ineq1} and \eqref{ineq2} with Theorem~\ref{thm:BS}, we obtain for each of the energies  $\F_p\in \{\E_p,\E_p^{\text{sym}},\I_p\}$, each $p\ge 2$ and each nontrivially knotted curve $\gamma$
\begin{equation}\label{ineq-chain-acn}
\F_p(\gamma)^{2/p}\ge \E_2^{\text{sym}}(\gamma)\ge 16\pi\cdot \text{acn}\,(\gamma) \ge 48\pi\, ,
\end{equation}
whereas for the circle    of length 1 (hence, radius $1/2\pi$) we have
\[
\E_p(\circle )^{2/p}=\E_p^{\text{sym}}(\circle )^{2/p}=\I_p (\circle )^{2/p}=(2\pi)^2 = 4\pi^2< 16\pi.
\]
The proof is complete now. \qed
\end{proof}

\begin{REMARK}
(i)\,
Instead of \eqref{ineq-chain-acn} we could have written
$$
\I_p(\gamma)^{2/p}\ge \I_2(\gamma)\ge\E_2(\gamma)\ge \E_2^{\text{sym}}(\gamma)\ge 16\pi\cdot \text{acn}\,(\gamma) ,
$$
and sending $p\to\infty$ does two things. Firstly,  it reproves one part of \cite[Theorem 4]{buck-simon_1997}, namely the inequality
\begin{equation}\label{thickness-e2}
\Big(\frac{1}{\triangle[\gamma]}\Big)^2\ge \E_2(\gamma).
\end{equation}
Secondly, it provides
 the lower ropelength bound   (as stated in \cite[Corollary 4.1]{buck-simon_1997})
\begin{equation}\label{thick-acn}
\Big(\frac{1}{\triangle[\gamma]}\Big)^2\ge 16\pi\cdot \text{acn}\,(\gamma)\ge 48\pi
\end{equation}
for nontrivial knots,
which is not quite as good as  the lower  bound  $1/
\triangle[\gamma]\ge 5\pi$ for any non-trivial knot obtained in
\cite[Corollary 3]{litherland-etal_1999}.

\noindent
(ii)\,
We can extend the inequality \eqref{rad_comparison} as
 $r_\textnormal{tp}[\gamma](x,y)\ge\varrho[\gamma](x,y)\ge\varrho_G[\gamma](x)\ge\triangle[\gamma],$
which implies the following order of energies for any unit loop $\gamma\in\cC$
(complementing the order in \eqref{energy-order} mentioned in the introduction):
$$
(\esym_p)^{1/p}(\gamma)\le\E_p^{1/p}(\gamma)\le\I_p^{1/p}(\gamma)\le\U_p^{1/p}(\gamma)\le\frac{1}{\triangle[\gamma]} \quad\Foa p\ge 2,
$$
and in particular again \eqref{thickness-e2}.
\end{REMARK}

We end this section by showing that for $p$ sufficiently  large, there is no non-trivial knot
minimizing integral Menger curvature $\M_p$, or $\E_p,$ $\esym_p$, or $\I_p.$
\begin{theorem}[\textbf{Trivial minimizers for multiple integral energies}]\label{thm:3.6}
There is a universal constant $p_0$ such that for all $p\ge p_0$ any minimizer of
$\M_p$, $\E_p$, $\esym_p$, or $\I_p$ is unknotted.
\end{theorem}
\proof
We restrict our proof to $\M_p$, analogous arguments work for the other energies as well.
We start with a general observation due to H\"older's inequality. If there is a curve
$\gamma\in\cC$ with $\M_p(\gamma)\le\M_p(\circle)$ for some $p>1,$
then the same inequality holds
true for any $q\in [1,p)$.

Assume that for all $n\in\N$ , $n\ge 4$, there exist $p_n>n,$ $p_{n+1}>p_n,$ and a non-trivially knotted
simple curve $\gamma_n\in\cC$ minimizing $\M_{p_n}$ in the class $\cC.$
Then in particular,
$$
\M_{p_n}^{1/p_n}(\gamma_n)\le\M_{p_n}^{1/p_n}(\circle)=2\pi\quad\Foa n\ge 4,
$$
so that we can use our initial  remark for $\gamma:=\gamma_n$, $p:=p_n>4$ and $q:=4$
to obtain
$$
\M_4^{1/4}(\gamma_n)\le 2\pi \quad\Foa n\ge 4.
$$
According to \cite[Theorem 4.3]{ssvdm-triple} this implies the
uniform a priori estimate
$$
\|\gamma_n\|_{C^{1,\alpha}([0,1],\R^3)}\le C\quad\Foa n\ge 4,
$$
where $\alpha=(4-3)/(4+6)=1/10.$
Hence there is a subsequence (still denoted by $\gamma_n$) converging in the $C^1$-norm
to a simple $C^1$-curve $\gamma_\infty\in\cC$ with finite energy $\M_4(\gamma)$, since $\M_p$
is lower-semicontinuous with respect to $C^1$-convergence (cf. Corollary 4.4 and Remark 4.5
in \cite{ssvdm-triple}).

{\it We claim that $\gamma_\infty$ is a circle of unit length.} Once this  is shown
we know by the isotopy result, Theorem \ref{thm:isotopy}, that $\gamma_n$ is unknotted for sufficiently large
$n$ contradicting our initial assumption, which proves the theorem.

Indeed,  we can estimate by lower semi-continuity of
$\M_{p_n}$ for arbitrary $n\ge 4$
\begin{eqnarray*}
\M_{p_n}^{1/p_n}(\gamma_\infty)
& \le & \liminf_{k\to\infty}\M_{p_n}^{1/p_n}(\gamma_k)\\
& \le & \liminf_{k\to\infty}\M_{p_n}^{1/p_n}(\circle)=\M_{p_n}^{1/p_n}(\circle),
\end{eqnarray*}
where we have used our initial remark for $\gamma:=\gamma_k$, $p:=p_k$ for $k>n$, and
$q:=p_n$ in the last inequality. Letting $n\to\infty$ and hence also $p_n\to\infty$
we find
$$
\frac{1}{\triangle[\gamma_\infty]}\le\frac{1}{\triangle[\circle]},
$$
which implies our claim since the circle uniquely minimizes ropelength.
\qed


\section{Isotopies to polygonal lines and crossing number bounds}
\label{sec:4}

In this section, we prove the following two results, alluded to in the introduction.
\begin{theorem}[\textbf{Finite energy curves and their polygonal models}]\label{thm:4.1}        Let $\gamma\in {\mathcal C}$ be simple and $0<E<\infty$. Assume one of the following:
    \begin{enumerate}
        \renewcommand{\labelenumi}{{\rm (\roman{enumi})}}
        \item $\M_p(\gamma)\le E$ for some $p>3$;
        \item $\F_p(\gamma)\le E$ for some $p>2$, where $\F_p\in \{\I_p,\E_p,\E_p^{\text{sym}}\}$;
        \item  $\U_p(\gamma)\le E$ for some $p>1$.
    \end{enumerate}
Then, there exist constants $\delta_1=\delta_1(p)\in (0,1)$ and  $\beta=\beta(p)>0$ such that  $\gamma$ is ambient isotopic to the polygonal line $\bigcup_{i=1}^N[x_i,x_{i+1}]$
for each choice of points
\[
x_i=\gamma(t_i), \quad 0=t_1<t_2< \ldots < t_N, \quad x_{N+1}=x_1, \]
that satisfy
\[
 |x_i-x_{i+1}|< \delta_1(p) E^{-\beta} \qquad\mbox{for $i=1,2,\ldots,N$}.
\]
We can take $\beta=1/(p-3)$ in case {\rm (i)},  $\beta=1/(p-2)$ in case {\rm (ii)}, and $\beta=1/(p-1)$ in case {\rm (iii)}.
\end{theorem}
 As an immediate consequence we note the following bound on the {\it stick number $\textnormal{seg}[K]$} of an isotopy class $[K]$,
 i.e.,
on the minimal number of segments needed to construct a polygonal
representative of $[K].$
\begin{corollary}[\textbf{Stick number}]\label{cor:stick_number}
Let $\gamma\in\cC$ be a representative of a knot class $[K]$, satisfying at least one
of the conditions (i), (ii), or (iii) in Theorem \ref{thm:4.1}. Then
\begin{equation}\label{stick_number_bound}
\textnormal{seg}[K]\le \frac{E^{\beta(p)}}{\delta_1(p)}+1.
\end{equation}
\end{corollary}
Since stick number and minimal crossing number are  strongly related (see, e.g.,\cite[Lemma 4]{litherland-etal_1999}) one immediately deduces an alternative
direct proof of the fact that all energies in Theorem \ref{thm:4.1} are strong
for the respective range of the parameter $p$,
and one could use the results in \cite[Section 3]{FHW} to produce explicit bounds
on the number of knot-types under a given energy level.
\begin{corollary}[\textbf{Finiteness}]\label{cor:finiteness}
Given $E>0$ and $p>1$, there can be at most finitely many knot types $[K_i]$ such that
there is a representative $\gamma_i\in\cC$ of $[K_i]$ with $\M_p(\gamma_i)\le E$
if  $p>3$, or with $\E_p(\gamma), \esym_p(\gamma), $ or $\I_p(\gamma)$
$\le E$ if  $p>2$, or with $\U_p(\gamma)\le E$.
\end{corollary}

\begin{theorem}[\textbf{Hausdorff distance related to energy implies isotopy}]\label{thm:4.2}      Let $\gamma_1,\gamma_2\in {\mathcal C}$ and $0<E<\infty$. Assume one of the following:
    \begin{enumerate}
        \renewcommand{\labelenumi}{{\rm (\roman{enumi})}}
        \item $\M_p(\gamma_j )\le E$ for some $p>3$ and $j =1,2$;
        \item $\F_p(\gamma_j )\le E$ for some $p>2$ and $j =1,2$, where $\F_p\in \{\I_p,\E_p,\E_p^{\text{sym}}\}$;
        \item  $\U_p(\gamma)\le E$ for some $p>1$.
    \end{enumerate}
Then, there exists a $\delta_2=\delta_2(p)\in (0,1)$  such that the two curves $\gamma_1$ and $\gamma_2$  are ambient isotopic if their Hausdorff distance does not exceed $\delta_2(p)E^{-\beta}$, with $\beta=1/(p-3)$ in case {\rm (i)}, $\beta=1/(p-2)$ in case {\rm (ii)},  and $\beta=1/(p-1)$ in case {\rm (iii)}.
\end{theorem}

For $x\not=y\in\bbbr^3$
and $\varphi\in (0,\frac\pi 2)$ we denote by
\[
C_\varphi(x;y)\, :=\, \{z\in\bbbr^3\setminus\{x\}\colon \exists \ t\neq 0  \  \
\textrm{such that}  \  \ \ang(t(z-x),y-x) < \frac \varphi 2 \} \cup\{x\}
\]
the double  cone whose  vertex is at the point $x$, with cone axis
passing through $y$, and with opening angle $\varphi$.

\begin{definition}[\textbf{Diamond property}]\label{def:diamonds}\label{def:4.0}
We say that a curve $\gamma\in {\mathcal C}$ has the \emph{diamond property at scale $d_0$ and with angle $\varphi\in (0,\pi/2) $}, in short \emph{the $(d_0,\varphi)$--diamond property}, if and only if for each couple of points $x,y\in \gamma$ with $|x-y|=d\le d_0$ two conditions are satisfied: we have
\begin{equation}
    \gamma\cap B_{2d}(x)\cap B_{2d}(y)\ \subset\ C_\varphi(x;y) \cap C_\varphi(y;x)
\end{equation}
(cf. Figure~\ref{2cones} below), and moreover each plane $a+(x-y)^\perp$, where $a\in B_{2d}(x)\cap B_{2d}(y)$, contains exactly one point of $\gamma\cap  B_{2d}(x)\cap B_{2d}(y)$.
\end{definition}

\begin{figure}[!h]

\begin{center}
\includegraphics*[totalheight=6.5cm]{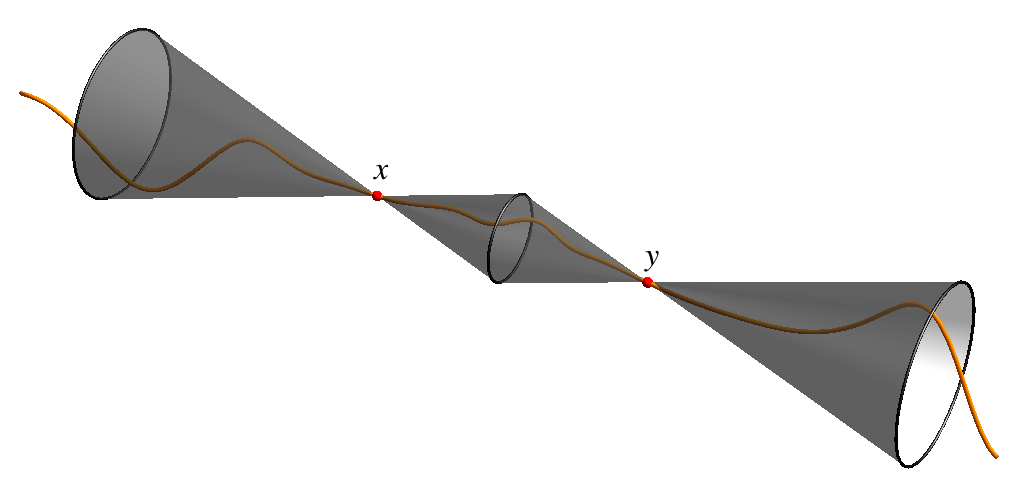}
\end{center}

\caption{\label{2cones}
The $(d_0,\varphi)$--diamond property: at small scales, the curve is trapped in a conical region and does not meander back and forth: each cross section of the cones contains exactly one point of the curve.   }
\end{figure}

Before proceeding further, let us note one immediate consequence of this property.

\begin{lemma}[\textbf{Bi-lipschitz estimate}] Suppose  a simple
curve $\gamma\in {\mathcal C}$ has  the $(d_0,\varphi)$--diamond property with $\varphi<1$. Then, whenever $|\gamma(s)-\gamma(t)|\le d_0\, $ for $|s-t|\le 1-|s-t|$, we have   \label{34-bilip}
    \[
    |\gamma (s)-\gamma(t)|\ge (1-\varphi) |s-t|\, .
    \]
\end{lemma}

\begin{proof} It is a simple argument, see e.g. \cite[Prop. 4.1]{svdm-tpcurves}. Assume first that $s\in [0,1]$ is a point of differentiability of $\gamma$. W.l.o.g suppose that $s<t$ and estimate
    \begin{eqnarray*}
    |\gamma(t)-\gamma(s)| & \ge &
    (\gamma(t)-\gamma(s))\cdot \gamma'(s)\\ & = &
    \int_s^t \bigl(\gamma'(\tau)-\gamma'(s)+\gamma'(s)\bigr)\,
    d\tau\cdot
    \gamma'(s) \\
    & \ge & (t-s) \Big(1 - \sup_{\tau\in [s,t]}
    |\gamma'(\tau)-\gamma'(s)|\Big)  \ \ge \ (1-\varphi) (t-s)\, .
    \end{eqnarray*}
(To verify the last inequality, let $S$ be the closed slab bounded
by two planes passing through $x=\gamma(s)$ and $y=\gamma(t)$, and perpendicular to $x-y$, i.e., to the common axis of the two cones; note that for each $\tau \in [s,t]$ we have in fact $\gamma(\tau)\in C_\varphi(x,y)
    \cap C_\varphi(y,x)\cap S$. Thus, for  all such $\tau$'s, we have $|\gamma'(s)-\gamma'(\tau)| \le\varphi$, as both vectors are of unit length and belong to the same  double cone with tips at $\gamma(s)$
    and $\gamma(\tau)$ and  opening angle  $\varphi$.)

Since the points of differentiability of $\gamma$ are dense in $[0,1]$, the lemma follows easily.
    \qed
\end{proof}

As we shall see, the diamond property allows to control the geometric behaviour (in particular, the bending at small and intermediate scales -- we will come to that later) of the curve. The main point is that finiteness of $\M_p$ (for $p>3$) or any one of the energies $\I_p$, $\E_p$ or $\E_p^{\text{sym}}$ (for $p>2$) implies the existence of two positive numbers $\alpha(p)$ and $\beta(p)$ such that each curve $\gamma\in{\mathcal C}$ of finite energy has the $(d_0,\varphi)$--diamond property at all sufficiently small scales $d_0\lesssim E^{-\beta}$ (where $E$ stands for the energy bound) with angle $\varphi \lesssim d_0^\alpha\ll 1$. Here is a more precise statement.

\begin{proposition}[\textbf{Energy bounds imply the diamond property}]\label{prop diam} Let $\gamma\in {\mathcal C}$ and $0<E<\infty$. Assume one of the following:
    \begin{enumerate}
        \renewcommand{\labelenumi}{{\rm (\roman{enumi})}}
        \item $\M_p(\gamma)\le E$ for some $p>3$;
        \item $\F_p(\gamma)\le E$ for some $p>2$, where $\F_p\in \{\I_p,\E_p,\E_p^{\text{sym}}\}$;
         \item  $\U_p(\gamma)\le E$ for some $p>1$.
    \end{enumerate}
Then, there exist constants $\delta=\delta(p)\in (0,1)$, $\alpha=\alpha(p)>0$, $\beta=\beta(p)>0$ and $c(p)<\infty$ (all four depending only on $p$) such that $\gamma$ has the $(d_0,\varphi)$--diamond property for each couple of numbers $(d_0,\varphi)$ satisfying
\begin{gather}
d_0\le  \delta(p) E^{-\beta}, \qquad \varphi \ge c(p) E^{\alpha\beta} d_0^\alpha\, .\label{d0 varphi}
\end{gather}
Specifically, we can take $\beta=1/(p-3)$, $\alpha=(p-3)/(p+6)$ in case {\rm (i)},  $\beta=1/(p-2)$, $\alpha=(p-2)/(p+4)$ in case {\rm (ii)}, and $\beta=1/(p-1)$, $\alpha=(p-1)/(p+2)$ in case {\rm (iii)}.
\end{proposition}

The proof of this proposition can be easily obtained from our earlier work (see \cite[Section 2]{ssvdm-triple} for the case of $\M_p$, \cite[Section~3]{ssvdm-double} for the case of $\I_p$, \cite[Section 4]{svdm-tpcurves} for the case of $\E_p$) and Kampschulte's master's thesis \cite{kampschulte_2012} for the case of $\E_p^{\text{sym}}$.  The last case of $\U_p$ can be treated via an application of \cite[Remark 7.2 and Theorem 7.3]{ssvdm-triple}, as the finiteness of $\U_p(\gamma)$ for $p>1$ and a simple curve $\gamma\in \cC$  implies, by H\"{o}lder inequality,
\[
\iiint_{B_r(\tau_1)\times B_r(\tau_2)\times B_r(\tau_3)} \frac{ds\, dt\, d\sigma}{R(\gamma(s),\gamma(t),\gamma(\tau))} \le 8r^{2+\delta} \U_p(\gamma), \qquad \delta= 1-\frac 1p >0,
\]
which is condition (7.2) of \cite{ssvdm-triple}.

In the remaining part of this section we will be working with double cones positioned along the curve. Let us introduce some notation first. For $x\not=y\in \R^3$ we denote the
closed halfspace
\begin{equation}
\label{bead}
H^+(x;y) \colon= \{z\in \R^3\colon \langle z-x, y-x\rangle \ge
0\}\, ,
\end{equation}
and use the `double cones'
\begin{equation}
    \label{double c}
    K(x,y)\colon= C_{1/4}(x;y) \cap C_{1/4}(y;x) \cap H^+(x;y) \cap
    H^+(y;x)\, .
\end{equation}

                     \begin{figure}[!t]

\begin{center}
\includegraphics*[totalheight=8.5cm]{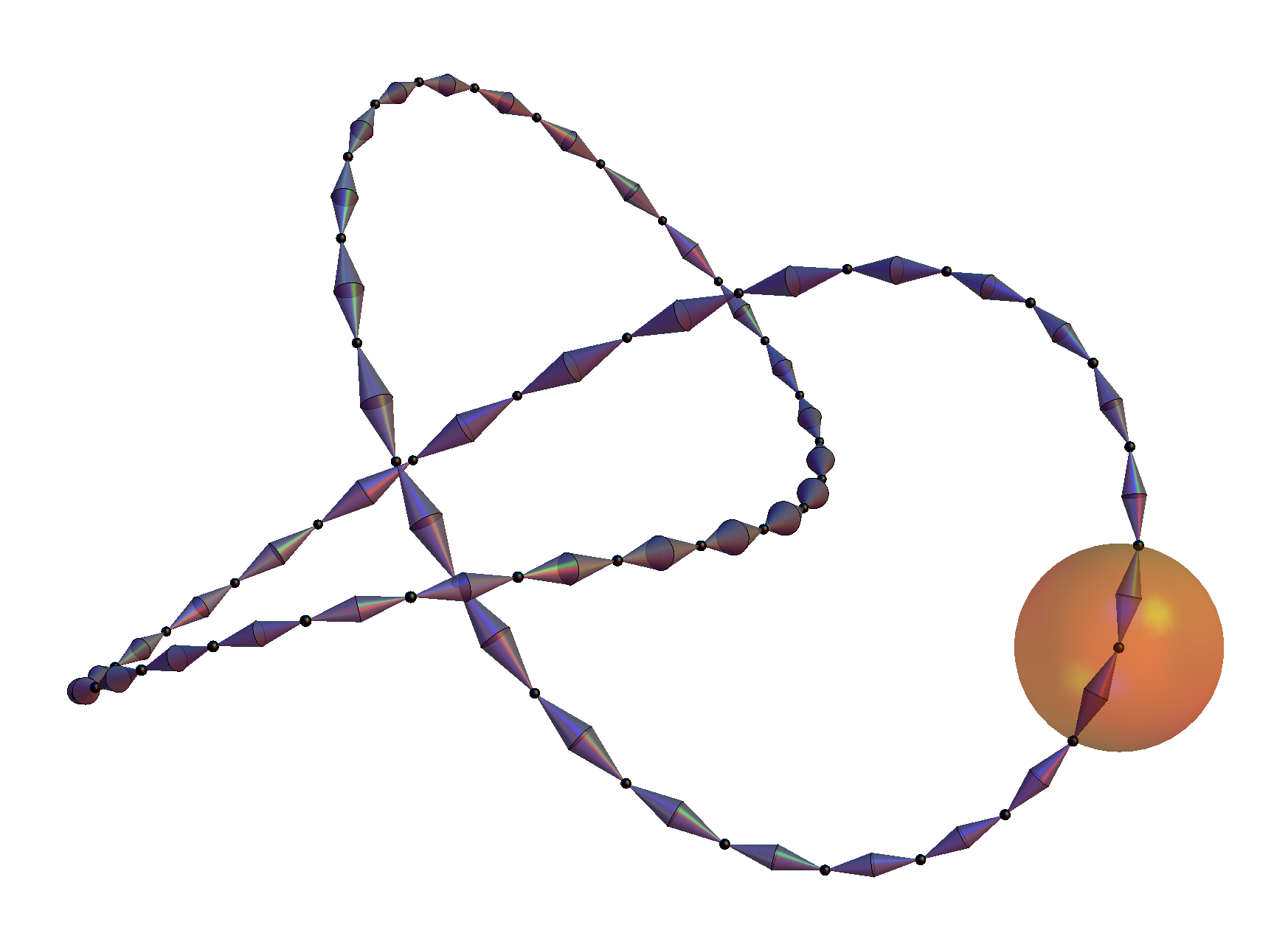}
\end{center}

\caption{\label{necklace1}
The meaning of the ``Necklace Lemma~\ref{lemma necklace}'': small double cones with vertices along the curve have pairwise disjoint interiors. Moreover, it follows from the $(d_0,\frac 14)$--diamond property that the angle between the axes of two neighbouring cones is at most $\frac 18$, and that different portions of the necklace stay well away from each other. If the points $x_i$ are evenly spaced, $|x_i-x_{i+1}|\equiv d_0$, then each ball $B_{d_0}(x_i)$ contains only the arcs of $\gamma$ coming from the two
double cones with common vertex at $x_i$.      The polygonal curve isotopic to $\gamma$, cf. Theorem~\ref{metathm:B}, joins the consecutive vertices of the cones.  }
\end{figure}
         \begin{lemma}[\textbf{Necklace of disjoint double cones}]\label{lemma necklace}  Suppose that $\gamma\in {\mathcal C}$ is simple
and has the $(d_0,\frac 14)$--diamond property. If  $0= t_1<\ldots < t_N< 1$ and $t_{N+1}=t_1$ and $x_i=\gamma(t_i)$ are such that $|x_{i+1}-x_i|\le d_0$, then the open double cones
\[
K_i=\text{int}\, K(x_i,x_{i+1}) \qquad\mbox{and} \qquad K_j=\text{int}\, K(x_j,x_{j+1})
\]
are disjoint whenever $i\not=j \pmod N$. Moreover, the vectors $v_i=x_{i+1}-x_i$ satisfy $\ang (v_{i+1},v_i)< 1/8$.
\end{lemma}

\begin{REMARK} The number $1/4$ in the lemma has been chosen just for the sake of simplicity, in favour  of simple arithmetics used now instead of more complicated computations in the theorems that follow. The result holds in fact for any angle $\varphi \le \frac 14$, with $\frac 18$ replaced by $\varphi/2$.
\end{REMARK}

\begin{proof}  By the $(d_0,\frac 14)$--diamond property, for each $z\in K_i$ the intersection of $\gamma$ and the two-dimensional disk
    \[
    D_i(z):= K_i\cap (z+v_i^\perp)
    \]
    contains precisely one point. Now, suppose to the contrary that
    \begin{equation}\label{A}
    K_i\cap K_j\not=\emptyset,
    \end{equation}
    and assume without loss of generality
    \begin{equation}\label{B}
    \diam K_j\le\diam K_i.
    \end{equation}
    If $x_j=\gamma(t_j)$ were contained in $K_i$ then
   either the disk $D_i(x_j)$ would contain two distinct curve
    points contradicting the second condition of the diamond property, or there would
be a parameter $\tau\in (t_i,t_{i+1})$
    such that $\gamma(\tau)=\gamma(t_j)$ although $\gamma$ is injective, a contradiction.
    The same reasoning can be applied to $x_{j+1}=\gamma(t_{j+1}),$ so that we conclude
    from  assumptions \eqref{A} and \eqref{B} that the two
    tips
    $x_j,$ $x_{j+1}$ of $K_j$ are contained in the set $Z_i$ defined as
    \begin{equation}\label{C}
    Z_i:=C_{\frac 14}(x_i;x_{i+1})\cap C_{\frac 14} (x_{i+1};x_i)\cap
    B_{2|v_i|}(x_i)\cap B_{2|v_i|}(x_{i+1})\; \setminus\; K_i,
    \end{equation}
    which is just the intersection of the two cones within the balls centered
    in $x_i$ and $x_{i+1}$ but without the open slab bounded by the two parallel
    planes $\partial H^+(x_i,x_{i+1})$ and
    $\partial H^+(x_{i+1},x_i)$.


\omitted{\tt\yy\yy\yy wouldn't this be shorter and simpler -- 8 lines instead of 12, and the logic is more direct:}

Since $\diam K_j=|v_j|\le \diam K_i=|v_i|$, we either have $\{x_i,x_{i+1}\}=\{x_j,x_{j+1}\}$ which
 in combination with the diamond property clearly contradicts the injectivity of $\gamma$,
or both points $x_j,x_{j+1}$ are in the same connected component of $Z_i$,
say in the one contained in  $\R^3\setminus H^+(x_{i+1},x_i)$. To fix the ideas, suppose that
$x_j$ is closer to the plane $\partial H^+(x_{i+1},x_i)$ than $x_{j+1}$ (or both points are equidistant from that plane). Then, the segment $[x_j,x_{j+1}]$ is contained in $H^+(x_j,x_j+v_i)$ so that all points of $K_j$ are contained outside the infinite half-cone
\[
S:=C_{\pi-\frac 14}(x_j;x_j+v_i) \cap H^+(x_j;x_j-v_i),
\]
which clearly contradicts \eqref{A} since, as it is easy to see, $K_i\subset S$.

The condition $\ang(v_i,v_{i+1})< \frac 18$ follows directly from the diamond property: without loss of generality, reversing the orientation of $\gamma$ if necessary, we may suppose that $|x_{i+2}-x_{i+1}|\le |x_{i+1}-x_i|=:d$. Then, $v_{i+1}\in C_{1/4}(0;v_i)$, and the inequality follows. \hfill $\Box$
\end{proof}

\begin{theorem}[\textbf{Isotopies to polygonal lines}]\label{metathm:B} Suppose that $\gamma\in {\mathcal C}$ is simple and has the $(d_0,\frac 14)$--diamond property. Then $\gamma$ is ambient isotopic to the polygonal curve
    \[
    P_\gamma = \bigcup_{i=1}^N [x_i,x_{i+1}]
    \]
    with $N$ vertices $x_i=\gamma(t_i)$, whenever the
    parameters $0= t_1<\ldots < t_N< 1$ and $t_{N+1}=t_1$ are chosen in $[0,1]$ so that
    \begin{equation}\label{spacing}
    |x_i-x_{i+1}|< d_0\, , \qquad i=1,\ldots, N.
    \end{equation}
\end{theorem}

\begin{proof} To construct the isotopy from $\gamma$ to a polygonal curve, we rely on Lemma~\ref{lemma necklace} and the diamond property. Cover $\gamma$ with a necklace of double cones $K(x_i,x_{i+1})$ that have pairwise disjoint interiors. The desired isotopy is constant off the union of $K(x_i,x_{i+1})$, and on each double cone it maps each two dimensional cross section  $D_i(z):= K_i\cap (z+v_i^\perp)$,
 where $z\in [x_i,x_{i+1}]$ and $v_i=x_{i+1}-x_i$, homeomorphically to itself, keeping the boundary of $D_i(z)$ fixed and moving the point $\gamma(s)\in D_i(z)$ along a straight segment until it hits the axis of the cone. \hfill $\Box$

\end{proof}

\begin{theorem}[\textbf{Isotopy by Hausdorff distance}]\label{metathm:C} Suppose that two simple curves $\gamma_1,\gamma_2\in {\mathcal C}$ are of class $C^1$ and have the $(d_0,\frac 14)$--diamond property. If their Hausdorff distance is smaller than $\eps=d_0/150$ then $\gamma_1$ and $\gamma_2$ are ambient isotopic.
\end{theorem}

\begin{REMARK}
As in Theorem 1.2 of \cite{svdm-tpcurves} it is actually not necessary to assume equal length
of $\gamma_1$ and $\gamma_2$.
\end{REMARK}

\begin{proof} Fix $\eta = \frac 13 d_0$ and pick $N>1/\eta\ge N-1 $, so that
 $t_i:=(i-1)\eta\in [0,1]$, $i=1,\ldots,N+1$, with the standard convention $t_{N+1}=t_1$
yield an equidistant partition of $[0,1]$. Assume now that $\dist_{H}(\gamma_1,\gamma_2)< \eps=d_0/150$. By
Theorem~\ref{metathm:B}, $\gamma_1$ is ambient isotopic to the
polygonal line
\[
P_{\gamma_1}:=\sum_{i=1}^N [x_i,x_{i+1}]\, ,
\]
where $x_i:=\gamma_1(t_i)$. Now, for $i=1,\ldots,N$ we set
$w_i:=\gamma_1'(t_i)$,
$\alpha_i:=\gamma_1\bigl([t_i,t_{i+1}]\bigr)\subset\gamma_1$, and
introduce the half-spaces $H^+_i:=H^+(x_i,x_i+w_i)$ and
$H_i^-:=\R^3\setminus H_i^+$, which are bounded by affine planes
$P_i:=x_i+w_i^\perp$.

The goal of the proof is to select points $y_i\in \gamma_2$ in each of the $P_i$ so that the polygonal line $P_{\gamma_2}$ with vertices at the $y_i$ would be isotopic both to $\gamma_2$ (via Theorem~\ref{metathm:B}) and to $P_{\gamma_1}$ (via an appropriate sequence of $\Delta$ and $\Delta^{-1}$ moves).

Throughout the whole proof, $|\sigma-\tau|$ etc. always refers to the \emph{intrinsic\/} distance of parameters on the circle of length 1.

\smallskip\noindent\emph{Step 1. Disjoint tubular regions around $P_{\gamma_1}$.} 
Consider the tubular regions
\[
T_i:= H_i^+ \cap H_{i+1}^- \cap B_{18\eps}(\alpha_i) .
\]
Their union contains $\gamma_1=\bigcup\alpha_i$; we clearly have
$T_i\cap T_{i+1}=\emptyset$ as $\alpha_{i+1}\subset H_{i+1}^+$.
In fact, we claim that \emph{$T_i\cap T_j=\emptyset$ whenever $|i-j|\ge 1$.} To see
this, we will use Lemma~\ref{34-bilip} to prove
\begin{equation}
\label{rigid} \inf\{|\gamma_1(\tau)-\gamma_1(\sigma)|\ \colon\
(\sigma,\tau)\in [0,1]\times [0,1],  \ |\sigma-\tau|\ge \eta\}
\ge \frac 34 \eta=\frac{3}{4}\cdot 50\eps.
\end{equation}

Before doing so, let us conclude from \eqref{rigid}: If there
existed a point $z\in T_i\cap T_j$ with $|i-j|>1$, we could find
$\sigma\in [t_i,t_{i+1})$ and $\tau \in [t_j,t_{j+1})$ such that
$|\gamma_1(\sigma)-\gamma_1(\tau)|\le 36 \eps<150\epsilon/4=3\eta/4$ by the triangle
inequality, a contradiction to \eqref{rigid}.

To verify \eqref{rigid}, notice that Lemma~\ref{34-bilip} applied to $\gamma_1$ implies
\begin{equation}\label{plus}
|\gamma_1(\tau)-\gamma_1(\sigma)|\ge\frac 34 |\tau-\sigma|\ge
\frac 34 \eta \quad\Foa\eta\le |\tau-\sigma|\le 3\eta.
\end{equation}
Now, since $\gamma$ is injective on $[0,1)$, the continuously differentiable function $g:[0,1]^2\to\R$ given by $g(s,t):=|\gamma_1(s)-\gamma_1(t)|^2$
attains a positive minimum $g_0>0$ on the compact set
$K_{3\eta}$, where we set $K_\rho:=[0,1]^2\setminus\{|s-t|<\rho\}$. Let $(s^*,t^*)\in K_{3\eta}$ be such that
$g(s,t)\ge g(s^*,t^*)=g_0$ for all $(s,t)\in K_{3\eta}.$
If $|s^*-t^*|=3\eta$ we can apply \eqref{plus} to find
$$
|\gamma_1(\tau)-\gamma_1(\sigma)|=\sqrt{g(\tau,\sigma)}
\ge\sqrt{g(s^*,t^*)}=|\gamma_1(s^*)-\gamma_1(t^*)|\overset{\eqref{plus}}{\ge}
\frac 34 \eta\quad\Foa (\tau,\sigma)\in K_{3\eta}.
$$
If, on the other hand, $|s^*-t^*|>3\eta$ then by minimality
$\nabla g(s^*,t^*)=0$, which implies that both tangents $\gamma_1'(s^*)$
and $\gamma_1'(t^*)$ are perpendicular to the segment $\gamma_1(s^*)-
\gamma_1(t^*).$ Thus the intersection
$$
\gamma_1([0,1])\cap B_{2\sqrt{g_0}}(\gamma_1(s^*))
\cap B_{2\sqrt{g_0}}(\gamma_1(t^*))
$$
cannot be contained in the intersection
$C_{1/4}(\gamma_1(s^*),\gamma_1(t^*))\cap
C_{1/4}(\gamma_1(t^*),\gamma_1(s^*)),
$
which according to the diamond property means that
$$|\gamma_1(s^*)-\gamma_1(t^*)|> d_0 = 3\eta,
$$
thereby establishing \eqref{rigid} also in this case.

\smallskip\noindent\emph{Step 2.} To choose a polygonal line that is ambient isotopic to $\gamma_2$, we prove the following: {\it for each $i=1,\ldots,N$ there is a
point
$
y_i\in P_i\cap\gamma_2\cap B_{2\eps}(x_i).
$
}

\smallskip
Without loss of generality we can assume that the curve $\gamma_1$
is oriented in such a way that
\begin{equation}\label{orientation}
\ang (\gamma_1'(t_i),v_i)<\frac 18 \quad\AND\quad\ang(\gamma_1'(t_i),v_{i-1})<\frac 18\quad\Foa i=1,\ldots,N,
\end{equation}
that is, each tangent $\gamma_1'(t_i)$ points into the set
$K_i:=K(x_i,x_{i+1})=K(\gamma_1(t_i),\gamma_1(t_{i+1})),$
which readily implies for the hyperplanes $P_i\perp\gamma_1'(t_i)$,
$i=1,\ldots,N$,
$$
\ang (P_i,v_i)\ge\ang (P_i,\gamma_1'(t_i))-\ang(\gamma_1'(t_i),v_i)
>\frac{\pi}{2}-\frac 18,
$$
and similarly $\ang(P_i,v_{i-1}) > \frac{\pi}{2}-\frac 18$. Indeed, according to the diamond property,
$$
\Big[\gamma_1\cap B_{2|v_i|}(x_i) \cap B_{2|v_i|}(x_{i+1})\cap
H^+(x_i,x_{i+1})\cap H^+(x_{i+1},x_i)\Big] \subset K_i,
$$
which implies that the tangent direction of the curve $\gamma_1$
at $x_i$ cannot deviate too much from the straight line through
$x_i$ and $x_{i+1}$; the inequalities in \eqref{orientation}
provide a quantified version of this fact.

Since $\dist_H(\gamma_1,\gamma_2)<\eps$ we find three points
$$
z_i\in\gamma_2\cap B_\eps(x_i),\quad
z_{i+1}\in\gamma_2\cap B_\eps(x_{i+1})\quad\AND\quad
z_{i-1}\in\gamma_2\cap B_\eps(x_{i-1})\quad\Foa i=1,\ldots,N.
$$
If $z_i\in P_i$ we set $y_i:=z_i$, and we are done. Else we know
that $z_i\in H^+_i\setminus P_i$ or that $z_i\in H_i^-.$
In the first case we will work with the two points $z_i$ and $z_{i-1}$,
in the second with $z_i$ and $z_{i+1}$ in the same way, so let us
assume the second situation $z_i\in H_i^-.$ We know that $z_{i+1}
\in  H_i^+\setminus P_i$ since by Lemma~\ref{34-bilip}
$$
\dist(z_{i+1},H^-_i)\ge\dist (x_{i+1},H_i^-)-\eps 
\ge
\left(\frac 34 - \frac 1{50}\right)\eta >0.
$$
On the other hand, $z_i$ and $z_{i+1}$ are not too far apart,
$$
\rho_i:=|z_i-z_{i+1}|\le|z_i-x_i|+|x_i-x_{i+1}|+|x_{i+1}-z_{i+1}|<
2\eps+\eta<d_0
$$
so that we can infer from the diamond property of $\gamma_2$ applied to the points
$x:=z_i$ and $y:=z_{i+1}$ that
\begin{equation}\label{PLUS}
\gamma_2\cap B_{2\rho_i}(z_i)
\cap B_{2\rho_i}(z_{i+1})\cap H^+(z_i,z_{i+1})\cap H^+(z_{i+1},z_i)
\subset K(z_i,z_{i+1}).
\end{equation}
We will now show that
\begin{equation}\label{PLUSPLUS}
\Big[ K(z_i,z_{i+1})\cap P_i\Big] \subset B_{2\eps}(x_i).
\end{equation}
Notice that $K(z_i,z_{i+1})\setminus P_i$ consists of two
components, one containing $z_i\in\gamma_2$, and the other
one containing $z_{i+1}\in\gamma_2,$ which implies
that the intersection in \eqref{PLUSPLUS} is not empty.
Since $\gamma_2$ connects $z_i$ and $z_{i+1}$
 by \eqref{PLUS} within the set $K(z_i,z_{i+1})$,
the inclusion in \eqref{PLUSPLUS} yields the desired curve point
$$
y_i\in P_i\cap\gamma_2\cap B_{2\eps}(x_i)\quad\Foa i=1,\ldots,N,
$$
thus proving the claim.

To prove \eqref{PLUSPLUS} we first estimate the angle $\ang(z_{i+1}-
z_i,v_i)$ by the largest possible angle between a line tangent
to both  $B_\eps(x_i)$ and $B_\eps(x_{i+1})$ and the line connecting
the centers $x_i,$ $x_{i+1}$:
$$
\ang(z_{i+1}-z_i,v_i)\le \arcsin\frac{\eps}{|v_i|/2},
$$
so that, using \eqref{orientation} and the estimate of  $|v_i|$ that follows from Lemma~\ref{34-bilip},
\[
\ang(z_{i+1}-z_i,\gamma_1'(t_i))  <  \frac 18
+\arcsin\frac{2\eps}{|v_i|}
 < \frac 18
+\arcsin\frac{2\eta/50}{3\eta/4}  <  \frac 15.
\]
Now, let $\tilde{z_i}$ be the orthogonal projection of $z_i$ onto
$P_i$. Since $\ang (\tilde{z}_i-z_i,z_{i+1}-z_i) = \ang
(\gamma'(t_i),z_{i+1}-z_i) < \frac 15$, it is easy to see that $
K(z_i,z_{i+1}) \cap P_i \subset B_{\tilde h}(\tilde{z}_i) \cap P_i $
where
$$
\tilde{h}\le|z_i-\tilde{z}_i| \tan\Big(\frac 15 + \frac 18\Big)<
\eps\tan\Big(\frac{8+5}{40}\Big)<\frac {\eps}2
$$
(see Figure~\ref{z cone} below), which establishes $K(z_i,z_{i+1})\cap
P_i\subset B_{2\eps}(x_i)$ and hence \eqref{PLUSPLUS}.

\begin{wrapfigure}[16]{l}[0cm]{6.9cm}

\includegraphics*[totalheight=4.1cm]{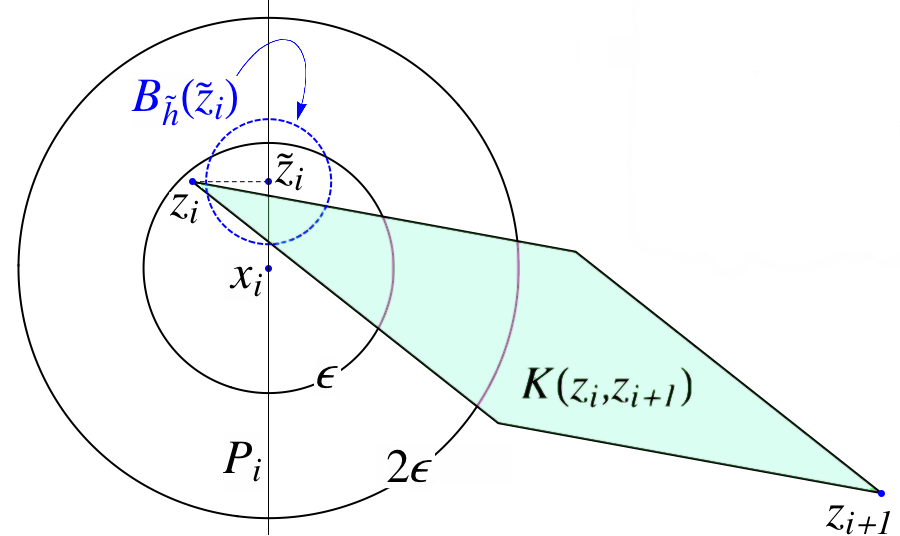}

\vspace*{1.5mm}

\caption{\label{z cone} The intersection
of the double cone $K(z_i,z_{i+1})$ with the plane $P_i$
is contained in $B_{\tilde h}(\tilde{z}_i)\subset
B_{2\eps}(x_i)$.
}

\end{wrapfigure}

Since $|y_i-y_{i+1}|< \eta + 4\eps< 3\eta =d_0$, the curve $\gamma_2$ is ambient isotopic to the
polygonal curve $P_{\gamma_2} =\bigcup_{i=1}^N [y_i,y_{i+1}]$.

\smallskip\noindent\emph{Step 3.} To
finish the proof of Theorem~\ref{metathm:C}, it is now sufficient to
check that {\it $P_{\gamma_1} $ and $P_{\gamma_2} $ are
combinatorially equivalent.} Since the sets $T_i$ are pairwise disjoint according
to Step~1, and
\[
B_{5\eps}\bigl([x_i,x_{i+1}]\bigr)\cap H^+_i\cap H^-_{i+1}\ \subset \ T_i,
\]
\heikodetail{\tt\xx this last inclusion should still be true even with the smaller choice of
width for $T_i$ by computing the height of $K_i\supset\alpha_i$ equal
$|(v_i|/2)\tan\frac 18 <\eta/14=50\epsilon/14$ as an upper bound for
$\dist (\alpha_i,[x_i,x_{i+1}])$.}
we have
\[
\mathrm{conv} (x_i,x_{i+1},y_i,y_{i+1}) \cap P_{\gamma_1} =
[x_i,x_{i+1}].
\]
This guarantees that all steps in the construction that follows
involve legitimate $\Delta$ and $\Delta^{-1}$-moves.
 (For the definition of these moves, and the distinction between them and the so-called Reidemeister moves, we refer to Burde and Zieschang's monograph \cite[Chapter~1]{burde-zieschang}). The first
step, taking place in $\overline{T}_1$, is to replace $[x_1,x_2]$
by the union of $[x_1,y_1]$ and $[y_1,x_2]$, and then to replace
$[y_1,x_2]$ by the union of $[y_1,y_2]$ and $[y_2,x_2]$. Next we
perform one $\Delta^{-1}$ and one $\Delta$-move in each of the
$\overline{T}_j$ for $j=2,\ldots, N-1$, replacing first
$[y_j,x_j]$ and $[x_{j},x_{j+1}]$ by $[y_j,x_{j+1}]$, and next
trading $[y_j,x_{j+1}]$ for the union of  $[y_j,y_{j+1}]$ and
$[y_{j+1},x_{j+1}]$. Finally, for $j=N$ we perform two
$\Delta^{-1}$-moves: first replace $[y_N,x_N]$ and $[x_{N},x_{1}]$
by $[y_N,x_1]$, and then replace $[y_N, x_1]$ and $[x_1,y_1]$
(which has been added at the very beginning of the construction)
by $[y_N,y_1]$. This concludes the whole proof. \qed

\end{proof}

{\sc Proof of Theorems \ref{thm:4.1} and \ref{thm:4.2}.}\,
 For $E$ fixed and $d_0\to 0$ condition \eqref{d0 varphi} of Proposition~\ref{prop diam} gives angles $\varphi\approx d_0^\alpha\to 0$. As we have already noted, this observation can be used to prove that all curves with bounded $\M_{p>3}$, $\I_{p>2}$, $\E_{p>2}$, or $\E_{p>2}^{\text{sym}}$ energy are in fact $C^1$, even $C^{1,\alpha}$ for some $\alpha>0$. Therefore, both Theorem~\ref{metathm:B} and Theorem~\ref{metathm:C} can be used for these energies; in combination with Proposition~\ref{prop diam} this clearly yields the two theorems stated at the beginning of this section.
    \qed

We end this section with a crude estimate of the average crossing number for curves that have the diamond property.

\begin{proposition}
\label{acn est}
Let $\gamma \in \cC$. Assume that there exists $d_1$ such that for each $d\le d_1$ the curve $\gamma$
satisfies the $\big(d,\varphi(d)\big)$-diamond property, where $\varphi(d) = Cd^\alpha$ for some $\alpha\in (\frac{1}{2},1]$ and $\varphi(d_1)\le \frac{1}{4}$ .
Then the average crossing number of the curve is finite and there exist absolute constants $\xi_1$ and $\xi_2$ such that
\begin{equation}
     \label{diam Acn}
    \mbox{{\em acn}}(\gamma) < \frac{C^2\xi_1 }{2\alpha-1} d_1^{2\alpha-1}  + \xi_2 d_1^{-\frac{4}{3}}
\end{equation}
\end{proposition}


\medskip

The general idea of the proof of Proposition \ref{acn est} is analogous to \cite[Cor.~4.1]{buck-simon_1997} and \cite[Cor.~2.1]{buck-simon_1999}.
We split the integral expressing the average crossing number into two parts; one of them, the local contribution, can be controlled using the local smoothness properties of the curve; the other one takes into account the interactions of distant portions of the curve. The novelty here is that the diamond property can be used to provide an excluded volume constraint and bound the length of the curve in a spherical shell around each of its points.

\medskip

\begin{proof}
First we notice that the expression in the numerator of the integrand of \eqref{acn 2} is equal to the volume of the parallelepiped spanned by vectors $\gamma'(s), \gamma'(t)$ and $\gamma(s) - \gamma(t)$. For a curve which satisfies the $(d, \varphi)$-diamond property the angles between the derivatives, and the derivatives and the secant, can be easily estimated.  Thus, for $|\gamma(s) - \gamma(t)| \le d_1$, we obtain, proceeding as in the proof of Theorem~\ref{thm:BS},
\[
|(\gamma'(s)\times\gamma'(t))\cdot (\gamma(t)-\gamma(s))|=
|\det(\gamma'(s), \gamma'(t), \gamma(s) - \gamma(t))| \le |\gamma'(s)| |\gamma'(t)| |\gamma(s) - \gamma(t)| \sin \varphi \sin \frac{\varphi}{2},
\]
where, by assumption, we can use $\varphi= C |\gamma(s) - \gamma(t)|^\alpha$. Hence,
\begin{equation}
\label{cross est}
\frac{|\det(\gamma'(s), \gamma'(t), \gamma(s) - \gamma(t))|}{|\gamma(s) - \gamma(t)|^3} \le
\frac{1}{2}C^2 |\gamma(s) - \gamma(t)|^{2\alpha - 2}.
\end{equation}

To estimate $\text{acn}(\gamma)$ we split the domain of integration into two parts.  We denote $S^1 = \R/\Z$ and  set
\[
X_s := \{t \in S^1 \ | \ |s-t|  \le  d_1 \}.
\]
Inequality \eqref{cross est} implies
\begin{align*}
I_X & :=  \int_{S^1} \int_{X_s} \frac{|\text{det}(\gamma'(s), \gamma'(t), \gamma(s) - \gamma(t))|}{|\gamma(s) - \gamma(t)|^3}dt\, ds  \\& \le \int_{S^1} \int_{X_s}\frac{1}{2}C^2 |\gamma(s) - \gamma(t)|^{2\alpha - 2} dt \,ds \\
& \le  \int_{S^1} \int_{s-d_1}^{s+d_1} \frac{1}{2}C^2\Big(\frac 34\Big)^{2\alpha - 2}|s-t|^{2\alpha - 2} dt\,ds \qquad\mbox{by Lemma~\ref{34-bilip}}\\
& \le  \frac{C^2}{2 \alpha - 1}\Big(\frac{3}{4}\Big)^{2\alpha - 2} d_1^{2\alpha - 1}
\le  \frac 43\cdot \frac{C^2}{2 \alpha - 1}   d_1^{2\alpha - 1}\qquad\mbox{ as $\alpha\in (\frac 12 , 1]$.}
\end{align*}
\heikodetail{
{\tt\xx Using $|s-t|\le d_1$ in the $t$-integration I seem to obtain
$C^2(3/4)^{2\alpha-2}d_1^{2\alpha-1}$ without actually integrating, does that look better, or am I making a stupid mistake?}

{\tt\yy\yy Stupid mistake: $2\alpha-2<0$ if $\alpha<1$, and the integrand is unbounded. This is the spot where we *need* $\alpha > \frac 12$.}
}
To estimate the integral on the remaining part of the domain $ S^1 \times Y_s$, where  $Y_s : =  S^1 \setminus X_s$, we notice that  for $t \in Y_s$  we have
\[
|\gamma(s) - \gamma(t) |  >  \frac{3}{4} d_1,
\]
for otherwise, according to Lemma \ref{34-bilip},  we would have $\frac 34 |s-t|\le |\gamma(s) - \gamma(t) |  \le \frac{3}{4} d_1$, a contradiction for $t\in Y_s$.
We define a family of sets, whose union contains $Y_s$:
\begin{align*}
Y_s^0 := & \{t \in Y_s \ | \  \frac{3}{4}d_1  < |\gamma(s) - \gamma(t)| \le d_1\}\, , \\
Y_s^n := & \{t \in Y_s \ | \  n d_1 < |\gamma(s) - \gamma(t)| \le (n+1) d_1 \} \quad \text{for } n\in \N.
\end{align*}
Since the length of $\gamma$ is finite, there exists $N= N(d_1)$ such that
\[
I_Y =  \int_{S^1} \int_{Y_s} \frac{|\det(\gamma'(s),\gamma'(t),\gamma(s)-\gamma(t))|}{|\gamma(s)-\gamma(t)|^3} dt \, ds\le \sum_{n=0}^{N}\int_{S_1}\int_{Y_s^n} |\gamma(s) - \gamma(t)|^{-2} dt \, ds.
\]
Now our aim is to estimate from above the measure of each $Y_s^n$.
We fix a polygonal curve with vertices
\[
x_i = \gamma(t_i) \Fo 0 = t_1 < t_2 < \ldots < t_k \quad\textnormal{with $x_{k+1} = x_{ 1}$},
\]
and with
\begin{equation}
\frac 34 d_1 \le |x_{i+1}-x_i| \le d_1\, .
\label{34d1}
\end{equation}
Since $\varphi(d_1)\le \frac 14$, $\gamma$ has the $(d_1,\frac 14)$--diamond property. Thus, by Lemma~\ref{lemma necklace},
\[
\gamma \subset \bigcup_{i\in I} K(x_i, x_{i+1}),  \qquad I=\{1,2,\ldots,k\},
\]
where $K(x,y)$ is the 'double cone' (with opening angle $\frac 14$) given by \eqref{double c}.
Using this inclusion we will find an upper bound for the length of the curve included in the spherical shells $A(a,b):= \overline{B_b(\gamma(s))}\setminus
\overline{B_a(\gamma(s))}$ for $0<a<b$  (if $a<0<b$ we simply put $A(a,b):= \overline{B_b(\gamma(s))}$).
\heikodetail{
{\tt\xx For consistency I changed notation for the
balls, but then the former sets $B_s$ needed to be changed, therefore I introduced
$Y_s$ (and  $X_s$) instead of $B_s$ (and $A_s$). Another issue, however, is the
correspondence of the spherical shells with the parameters described in the
former sets $B_s^n$ (not $Y_s^n$), we seem to need closed balls in the definition
of the shells, but this again does not seem to fit to $Y_s^0$. So, I gave up late at night, left it with the closure of the open balls here (balls are open in our paper so far), and leave the decision for you...\xx\xx}
}%
For fixed $a,b$, let $J \subseteq I$ denote
the set of all indices
$i \in I$ for which
\[
[\gamma\cap A(a,b)] \cap K(x_i, x_{i+1}) \neq \emptyset.
\]
Then we have
\[
\gamma\cap A(a,b) \subset \bigcup_{i \in J} K(x_i, x_{i+1}) \subseteq A(a-d_1, b+d_1).
\]
 Thus, the length of the portion of $\gamma$ within the spherical shell, measured
in the one-dimensional Hausdorff measure,  satisfies
\begin{equation}
    \label{H1 43}
    \H^1 \big(A(a,b) \cap \gamma\big) \le \H^1\big(\gamma \cap \bigcup_{i \in J}K(x_i, x_{i+1})\big) \le
    \sum_{i\in J}\frac{4}{3}|x_i- x_{i+1}|,
\end{equation}
where in the last inequality the bi-lipschitz continuity of the parametrization is used.

By Lemma \ref{lemma necklace} we know that $\text{int}\, K(x_i, x_{i+1}) \cap \text{int}\, K(x_j, x_{j+1}) = \emptyset$ for $i \neq j$. Thus we can estimate the volume of the union of 'double cones' from below
\begin{equation}
\label{vol below}
\H^{3}\big(\bigcup_{i \in J}K(x_i, x_{i+1})\big) =
\frac{\pi}{12} \tan^2{\frac{1}{8}}\sum_{i \in J}|x_i - x_{i+1}|^3 \stackrel{\eqref{34d1}}\ge
 \frac{\pi}{4^4}\Big(\frac{3 d_1}{4} \Big)^2  \sum_{i \in J} \frac{4}{3} |x_i - x_{i+1}|.
\end{equation}

On the other hand, the volume of $\bigcup_{i \in J}K(x_{i}, x_{i+1})$ cannot exceed the volume of $A(a-d_1, b+d_1)$. Therefore for $ a > d_1 $, combining \eqref{H1 43} and \eqref{vol below}, we obtain
\begin{align*}
\H^1 \big(A(a,b) \cap \gamma\big) \le \frac{4^4}{\pi} \Big(\frac{4}{3 d_1} \Big)^{2} \H^{3}\Big(\bigcup_{i \in J}K(x_1, x_{i+1}) \Big) \le 4^4\Big(\frac{4}{3}\Big)^3 d_1^{-2} [(b+d_1)^3  - (a-d_1)^3].
\end{align*}
Since $Y_s^n$ is just the preimage of $\gamma\cap A(nd_1,(n+1)d_1)$  for $n\ge 1$,  and $\gamma$ is simple and parametrized by arclength,
\begin{equation}
    \label{meas Bns}
    \H^1\big( Y_s^n \big) \le \frac{4^{7}}{3^3} [(n+2)^3 - (n-1)^3] d_1.
\end{equation}
Analogously, for $a < d_1$,
\[
\H^1 \big(A(a,b) \cap \gamma\big) \le \frac{4^{7}}{3^3}d_1^{-2} (b+d_1)^3,
\]
and (inserting $b:=d_1$)
\begin{equation}
    \label{meas B0}
    \H^1\big( Y_s^0\big) \le \frac{4^{7}}{3^3}  8 d_1.
\end{equation}

To estimate the integral $I_Y$ we assume the worst case which occurs when the curve is densely packed around the single point $\gamma(s)$ i.e. each shell $A(nd_1,(n+1)d_1)$ contains the maximum possible amount of length of the curve which is controlled by \eqref{meas Bns}. In this case, we can give an upper estimate  for $N=N(d_1)$, taking the smallest $N$ such that
\[
\sum_{n=0}^N \H^1(Y_s^n) \ge \H^{1}(\gamma) =1.
\]
Using \eqref{meas Bns}--\eqref{meas B0} we obtain
\begin{align*}
1 & \le\frac{4^7}{3^3} 8 d_1 +
 \sum_{n=1}^N \frac{4^7}{3^3} [(n+2)^3 - (n-1)^3] d_1
  = d_1 \frac{4^7}{3^3}\big((N+2)^3+(N+1)^3+ N^3 -1 \big).
\end{align*}
Thus it is enough to take the smallest integer $N$ such that $N^3 > \frac{3^2}{ 4^7}{d_1^{-1}}$. This gives the following estimation of the integral $I_Y$:
\begin{align*}
I_Y & \le \int_{S^1}\int_{Y_s^0} \Big(\frac{4}{3}\Big)^2 d_1^{-2} dt \, ds + \sum_{n=1}^N \int_{S^1}\int_{Y_s^n}(nd_1)^{-2} dt \, ds \\
&  \le d_1^{-1}  \Big(8 \frac{4^9}{3^5}  + \frac{4^7}{3^3} \sum_{n=1}^N  \frac{1}{n^2}[(n+2)^3 - (n-1)^3] \Big) \\
&\le 8d_1^{-1}\frac{4^9}{3^5} + \frac{4^7}{3^3} d_{1}^{-1} \cdot  3^3 N
\le  2d_1^{-1}\frac{4^{10}}{3^5} + \bigg(\Big(\frac{3^2}{ 4^7}{d_1^{-1}}\Big)^{1/3}+1\bigg) {4^7} d_1^{-1} < \xi d_1^{-\frac{4}{3}},
 \end{align*}
for some absolute constant $\xi$.

Eventually, we get the desired estimation for the average crossing number
\[
\text{acn}(\gamma) = \frac{1}{4\pi} I_X + \frac{1}{4\pi} I_Y \le \frac{C^2\xi_1 }{2\alpha-1} d_1^{2\alpha-1}+ \xi_2 d_1^{-\frac{4}{3}} \, .
\]
\qed
\end{proof}

Using Proposition \ref{prop diam} we get an estimate for the
average crossing number for the curves with finite $\M_p$ energy.

\begin{corollary} \label{our almost 43}
Let $\gamma \in \cC$ and $ 0<E<\infty$. If $\M_p(\gamma) < E$ for some $p>12$ then there exist constants $c_1(p)$ and $c_2(p)$, such that
\[
\text{\em{acn}}(\gamma) < c_1(p)
+ c_2(p)E^{\frac{4}{3(p-3)}}.
\]
\end{corollary}

\begin{proof}
According to Proposition \ref{prop diam},  we can express the constants $d_1$ and $C$ from Proposition \ref{acn est}  as
\[
d_1 =  \delta(p) E^{-\beta}, \qquad C =  c(p) E^{\alpha\beta}   ,
\]
where $\beta=1/(p-3)$ and $\alpha=(p-3)/(p+6)$. To obtain the required estimates, we insert  the above quantities into formula \eqref{diam Acn}, and next use the inequality $E^\beta\le 1+E^{4\beta/3}$.
\qed
\end{proof}


\omitted{\tt\xx just to make sure: just inserting all the quantities, in particular evaluating
the term $d_1^{2\alpha-1}$ would give some term $c(p)E^{\frac{p-12}{(p-3)(p+6)}}$, correct? \yy Not quite, as $C=C(E)$ in the Corollary, see the following text that Marta inserted.

 \omitted{$\heartsuit$} \omitted{$\heartsuit$}  In Proposition \ref{acn est} we assume that $\varphi(d) = Cd^\alpha$ and in estimation
\eqref{diam Acn} the first term depends on $C$ - namely  it is $\xi (\alpha) C^2 d_1^{2\alpha-1}$. Both $C$ and $d_1$ taken from  Proposition \ref{prop diam} depend on $p$  - inserting those values we get the exponent $\beta=1/(p-3)$ in the first term. \omitted{$\heartsuit$}  \omitted{$\heartsuit$}

Any simplification should also take into account
that $E$ might be quite small...
\xx}

\omitted{\yy\tt No power of $E$ is forgotten; as Marta noted, $C$ in (4.36) depends also on $E$ and since $\alpha$ depends only on $p$ this gives
\[
C^2 d_1^{2\alpha-1} = \text{const}(p) E^{2\alpha\beta} E^{-\beta(2\alpha-1)}=   \text{const}(p) E^\beta\le \text{const}(p) (1 + E^{4\beta/3})\, .
\]
}

\begin{REMARK} Since $\M_p(\gamma)^{1/p}$ approaches $1/\triangle[\gamma]$ as $p\to\infty$, and the constants $c_1(p)$, $c_2(p)$ do not blow up\footnote{This can be checked by tracing the constants in \cite{ssvdm-triple}.} as $p\to \infty$, Corollary~\ref{our almost 43} gives, in the limit $p\to \infty$, a result which qualitatively agrees with Buck and Simon's \cite[Cor.~4.1]{buck-simon_1997} estimate of the average crossing number by a constant multiple of $(1/\triangle[\gamma])^{4/3}$. Our constant $c_2(p)$ is (far) worse, though.
\end{REMARK}

\addcontentsline{toc}{section}{References}

{

\bibliography{menger}{}

\begin{thebibliography}{10}

\bibitem{abrams-etal_2003}
Aaron Abrams, Jason Cantarella, Joseph H.~G. Fu, Mohammad Ghomi, and Ralph
  Howard.
\newblock Circles minimize most knot energies.
\newblock {\em Topology}, 42(2):381--394, 2003.

\bibitem{auckly-sadun_1997}
David Auckly and Lorenzo Sadun.
\newblock A family of {M}\"obius invariant {$2$}-knot energies.
\newblock In {\em Geometric topology ({A}thens, {GA}, 1993)}, volume~2 of {\em
  AMS/IP Stud. Adv. Math.}, pages 235--258. Amer. Math. Soc., Providence, RI,
  1997.

\bibitem{banavar}
Jayanth~R. Banavar, Oscar Gonzalez, John~H. Maddocks, and Amos Maritan.
\newblock Self-interactions of strands and sheets.
\newblock {\em J. Statist. Phys.}, 110(1-2):35--50, 2003.

\bibitem{blatt-reiter_2012}
S.~{Blatt} and P.~{Reiter}.
\newblock Regularity theory for tangent-point energies: {T}he non-degenerate
  sub-critical case, 2012.
\newblock Preprint.

\bibitem{blatt-reiter-schikorra_2012}
S.~{Blatt}, P.~{Reiter}, and A.~{Schikorra}.
\newblock Hard analysis meets critical knots ({S}tationary points of the
  {M}oebius energy are smooth), 2012.
\newblock arXiv:1202.5426v2.

\bibitem{blatt_isop_2009}
Simon Blatt.
\newblock Note on continuously differentiable isotopies, 2009.
\newblock Preprint Nr. 34, Institut f{\"u}r Mathematik, RWTH Aachen University.

\bibitem{blatt-tp}
Simon Blatt.
\newblock The energy spaces of the tangent--point energies, 2011.
\newblock {P}reprint.

\bibitem{blatt-menger_2011}
Simon Blatt.
\newblock A note on integral {M}enger curvature, 2011.
\newblock {P}reprint.

\bibitem{blatt-k}
Simon Blatt and S{\l}awomir Kolasi{\'n}ski.
\newblock Sharp boundedness and regularizing effects of the integral {M}enger
  curvature for submanifolds.
\newblock {\em Adv. Math.}, 230(3):839--852, 2012.

\bibitem{blatt-reiter}
Simon Blatt and Philipp Reiter.
\newblock Does finite knot energy lead to differentiability?
\newblock {\em J. Knot Theory Ramifications}, 17(10):1281--1310, 2008.

\bibitem{blumenthal-menger_1970}
Leonard~M. Blumenthal and Karl Menger.
\newblock {\em Studies in geometry}.
\newblock W. H. Freeman and Co., San Francisco, Calif., 1970.

\bibitem{buck-simon_1997}
Gregory Buck and Jonathan Simon.
\newblock Energy and length of knots.
\newblock In {\em Lectures at {KNOTS} '96 ({T}okyo)}, volume~15 of {\em Ser.
  Knots Everything}, pages 219--234. World Sci. Publ., River Edge, NJ, 1997.

\bibitem{buck-simon_1999}
Gregory Buck and Jonathan Simon.
\newblock Thickness and crossing number of knots.
\newblock {\em Topology Appl.}, 91(3):245--257, 1999.

\bibitem{burde-zieschang}
Gerhard Burde and Heiner Zieschang.
\newblock {\em Knots}, volume~5 of {\em de Gruyter Studies in Mathematics}.
\newblock Walter de Gruyter \& Co., Berlin, second edition, 2003.

\bibitem{CFKSW06}
Jason Cantarella, Joseph H.~G. Fu, Rob Kusner, John~M. Sullivan, and Nancy~C.
  Wrinkle.
\newblock Criticality for the {G}ehring link problem.
\newblock {\em Geom. Topol.}, 10:2055--2116 (electronic), 2006.

\bibitem{CKS02}
Jason Cantarella, Robert~B. Kusner, and John~M. Sullivan.
\newblock On the minimum ropelength of knots and links.
\newblock {\em Invent. Math.}, 150(2):257--286, 2002.

\bibitem{david-survey}
Guy David.
\newblock Analytic capacity, {C}alder\'on-{Z}ygmund operators, and
  rectifiability.
\newblock {\em Publ. Mat.}, 43(1):3--25, 1999.

\bibitem{dudziak_2010}
James~J. Dudziak.
\newblock {\em Vitushkin's conjecture for removable sets}.
\newblock Universitext. Springer, New York, 2010.

\bibitem{EG}
Lawrence~C. Evans and Ronald~F. Gariepy.
\newblock {\em Measure theory and fine properties of functions}.
\newblock Studies in Advanced Mathematics. CRC Press, Boca Raton, FL, 1992.

\bibitem{fary_1949}
Istv{\'a}n F{\'a}ry.
\newblock Sur la courbure totale d'une courbe gauche faisant un n\oe ud.
\newblock {\em Bull. Soc. Math. France}, 77:128--138, 1949.

\bibitem{federer_1959}
Herbert Federer.
\newblock Curvature measures.
\newblock {\em Trans. Amer. Math. Soc.}, 93:418--491, 1959.

\bibitem{FHW}
Michael~H. Freedman, Zheng-Xu He, and Zhenghan Wang.
\newblock M\"obius energy of knots and unknots.
\newblock {\em Ann. of Math. (2)}, 139(1):1--50, 1994.

\bibitem{fukuhara_1988}
Shinji Fukuhara.
\newblock Energy of a knot.
\newblock In {\em A f\^ete of topology}, pages 443--451. Academic Press,
  Boston, MA, 1988.

\bibitem{fuller-vemuri_2010}
E.J. Fuller and M.K. Vemuri.
\newblock The {B}rylinski beta function of a surface, 2010.
\newblock arXiv:1012.4096v1 [math.DG].

\bibitem{gonzalez-delallave_2003}
O.~Gonzalez and R.~de~la Llave.
\newblock Existence of ideal knots.
\newblock {\em J. Knot Theory Ramifications}, 12(1):123--133, 2003.

\bibitem{gonzalez-maddocks-smutny}
O.~Gonzalez, J.~H. Maddocks, and J.~Smutny.
\newblock Curves, circles, and spheres.
\newblock In {\em Physical knots: knotting, linking, and folding geometric
  objects in {$\Bbb R^3$} ({L}as {V}egas, {NV}, 2001)}, volume 304 of {\em
  Contemp. Math.}, pages 195--215. Amer. Math. Soc., Providence, RI, 2002.

\bibitem{GM}
Oscar Gonzalez and John~H. Maddocks.
\newblock Global curvature, thickness, and the ideal shapes of knots.
\newblock {\em Proc. Natl. Acad. Sci. USA}, 96(9):4769--4773 (electronic),
  1999.

\bibitem{GMSvdM}
Oscar Gonzalez, John~H. Maddocks, Friedemann Schuricht, and Heiko von~der
  Mosel.
\newblock Global curvature and self-contact of nonlinearly elastic curves and
  rods.
\newblock {\em Calc. Var. Partial Differential Equations}, 14(1):29--68, 2002.

\bibitem{He}
Zheng-Xu He.
\newblock The {E}uler-{L}agrange equation and heat flow for the {M}\"obius
  energy.
\newblock {\em Comm. Pure Appl. Math.}, 53(4):399--431, 2000.

\bibitem{hermes-phd}
Tobias Hermes.
\newblock {\em Analysis of the first variation and a numerical gradient flow
  for integral {M}enger curvature}.
\newblock PhD thesis, RWTH Aachen University, 2012.
\newblock Available at
  http://darwin.bth.rwth-aachen.de/opus3/volltexte/2012/4186/.

\bibitem{hirsch}
Morris~W. Hirsch.
\newblock {\em Differential topology}.
\newblock Springer-Verlag, New York, 1976.
\newblock Graduate Texts in Mathematics, No. 33.

\bibitem{kampschulte_2012}
Malte~Laurens Kampschulte.
\newblock The symmetrized tangent-point energy.
\newblock Master's thesis, RWTH Aachen University, 2012.

\bibitem{skol2}
S{\l}awomir Kolasi{\'n}ski.
\newblock {\em Integral {M}enger curvature for sets of arbitrary dimension and
  codimension}.
\newblock PhD thesis, Institute of Mathematics, University of Warsaw, 2011.
\newblock arXiv:1011.2008v4.

\bibitem{kolasinski_2012}
S{\l}awomir Kolasi{\'n}ski.
\newblock {G}eometric {S}obolev-like embedding using high-dimensional
  {M}enger-like curvature, 2012, submitted.
\newblock arXiv:1205.4112v1; Trans. Amer. Math. Soc., accepted for publication.

\bibitem{ksvdm-w2p}
S{\l}awomir Kolasi{\'n}ski, Pawe{\l} Strzelecki, and Heiko von~der Mosel.
\newblock Characterizing ${W}^{2,p}$-submanifolds by $p$-integrability of
  global curvatures, 2012.
\newblock arXiv:1203.4688v2 [math.CA]; Geometric and Functional Analysis,
  accepted for puublication.

\bibitem{ksvdm-tp2}
S{\l}awomir Kolasi{\'n}ski, Pawe{\l} Strzelecki, and Heiko von~der Mosel.
\newblock Tangent-point repulsive potentials for a class of non-smooth
  $m$-dimensional sets in $\mathbb{R}^n$. {P}art {II}: Compactness and
  finiteness results, 2012.
\newblock In preparation.

\bibitem{kusner-sullivan-moebius_1998}
R.~B. Kusner and J.~M. Sullivan.
\newblock M\"obius-invariant knot energies.
\newblock In {\em Ideal knots}, volume~19 of {\em Ser. Knots Everything}, pages
  315--352. World Sci. Publ., River Edge, NJ, 1998.

\bibitem{leger}
Jean-Christophe L{\'e}ger.
\newblock {M}enger curvature and rectifiability.
\newblock {\em Ann. of Math. (2)}, 149(3):831--869, 1999.

\bibitem{mattila-lin_2001}
Yong Lin and Pertti Mattila.
\newblock {M}enger curvature and {$C^1$} regularity of fractals.
\newblock {\em Proc. Amer. Math. Soc.}, 129(6):1755--1762 (electronic), 2001.

\bibitem{litherland-etal_1999}
R.~A. Litherland, J.~Simon, O.~Durumeric, and E.~Rawdon.
\newblock Thickness of knots.
\newblock {\em Topology Appl.}, 91(3):233--244, 1999.

\bibitem{mattila-melnikov-verdera}
Pertti Mattila, Mark~S. Melnikov, and Joan Verdera.
\newblock The {C}auchy integral, analytic capacity, and uniform rectifiability.
\newblock {\em Ann. of Math. (2)}, 144(1):127--136, 1996.

\bibitem{mcatee_2004}
Jenelle~Marie McAtee.
\newblock Knots of constant curvature, 2004.
\newblock arXiv:math/0403089v1 [math.GT].

\bibitem{melnikov-verdera}
Mark~S. Melnikov and Joan Verdera.
\newblock A geometric proof of the {$L^2$} boundedness of the {C}auchy integral
  on {L}ipschitz graphs.
\newblock {\em Internat. Math. Res. Notices}, (7):325--331, 1995.

\bibitem{menger_1930}
Karl {M}enger.
\newblock Untersuchungen \"uber allgemeine {M}etrik.
\newblock {\em Math. Ann.}, 103(1):466--501, 1930.

\bibitem{milnor_1950}
J.~W. Milnor.
\newblock On the total curvature of knots.
\newblock {\em Ann. of Math. (2)}, 52:248--257, 1950.

\bibitem{nitsche_1971}
J.~C.~C. Nitsche.
\newblock The smallest sphere containing a rectifiable curve.
\newblock {\em Amer. Math. Monthly}, 78:881--882, 1971.

\bibitem{ohara_1991}
Jun O'Hara.
\newblock Energy of a knot.
\newblock {\em Topology}, 30(2):241--247, 1991.

\bibitem{ohara_1992b}
Jun O'Hara.
\newblock Family of energy functionals of knots.
\newblock {\em Topology Appl.}, 48(2):147--161, 1992.

\bibitem{ohara_1994}
Jun O'Hara.
\newblock Energy functionals of knots. {II}.
\newblock {\em Topology Appl.}, 56(1):45--61, 1994.

\bibitem{ohara-book_2003}
Jun O'Hara.
\newblock {\em Energy of knots and conformal geometry}, volume~33 of {\em
  Series on Knots and Everything}.
\newblock World Scientific Publishing Co. Inc., River Edge, NJ, 2003.

\bibitem{reiter_isop_2005}
Philipp Reiter.
\newblock All curves in a ${C}^1$-neighbourhood of a given embedded curve are
  isotopic, 2005.
\newblock Preprint Nr. 4, Institut f{\"u}r Mathematik, RWTH Aachen University.

\bibitem{scholtes_2011}
Sebastian Scholtes.
\newblock {F}or which positive $p$ is the integral {M}enger curvature
  $\mathcal{M}_p$ finite for all simple polygons?, 2011.
\newblock arXiv:1202.0504v1.

\bibitem{scholtes_2012}
Sebastian Scholtes.
\newblock Tangency properties of sets with finite geometric curvature energies.
\newblock {\em Fund. Math.}, 218(2):165--191, 2012.

\bibitem{heiko1}
Friedemann Schuricht and Heiko von~der Mosel.
\newblock Global curvature for rectifiable loops.
\newblock {\em Math. Z.}, 243(1):37--77, 2003.

\bibitem{heiko3}
Friedemann Schuricht and Heiko von~der Mosel.
\newblock Characterization of ideal knots.
\newblock {\em Calc. Var. Partial Differential Equations}, 19(3):281--305,
  2004.

\bibitem{simon_1996}
Jonathan Simon.
\newblock Energy functions for knots: beginning to predict physical behavior.
\newblock In {\em Mathematical approaches to biomolecular structure and
  dynamics ({M}inneapolis, {MN}, 1994)}, volume~82 of {\em IMA Vol. Math.
  Appl.}, pages 39--58. Springer, New York, 1996.

\bibitem{ssvdm-double}
Pawe{\l} Strzelecki, Marta Szuma{\'n}ska, and Heiko von~der Mosel.
\newblock A geometric curvature double integral of {M}enger type for space
  curves.
\newblock {\em Ann. Acad. Sci. Fenn. Math.}, 34(1):195--214, 2009.

\bibitem{ssvdm-triple}
Pawe{\l} Strzelecki, Marta Szuma{\'n}ska, and Heiko von~der Mosel.
\newblock Regularizing and self-avoidance effects of integral {M}enger
  curvature.
\newblock {\em Ann. Sc. Norm. Super. Pisa Cl. Sci. (5)}, 9(1):145--187, 2010.

\bibitem{StvdM1}
Pawe{\l} Strzelecki and Heiko von~der Mosel.
\newblock On a mathematical model for thick surfaces.
\newblock In {\em Physical and numerical models in knot theory}, volume~36 of
  {\em Ser. Knots Everything}, pages 547--564. World Sci. Publ., Singapore,
  2005.

\bibitem{StvdM2}
Pawe{\l} Strzelecki and Heiko von~der Mosel.
\newblock Global curvature for surfaces and area minimization under a thickness
  constraint.
\newblock {\em Calc. Var. Partial Differential Equations}, 25(4):431--467,
  2006.

\bibitem{svdm-MathZ}
Pawe{\l} Strzelecki and Heiko von~der Mosel.
\newblock On rectifiable curves with ${L}^p$-bounds on global curvature:
  Self-avoidance, regularity, and minimizing knots.
\newblock {\em Math. Z.}, 257:107--130, 2007.

\bibitem{svdm-surfaces}
Pawe{\l} Strzelecki and Heiko von~der Mosel.
\newblock Integral {M}enger curvature for surfaces.
\newblock {\em Adv. Math.}, 226:2233--2304, 2011.

\bibitem{svdm-tp1}
Pawe{\l} Strzelecki and Heiko von~der Mosel.
\newblock Tangent-point repulsive potentials for a class of non-smooth
  $m$-dimensional sets in $\mathbb{R}^n$. {P}art {I}: Smoothing and
  self-avoidance effects, 2011.
\newblock arXiv:1102.3642; J. Geom. Anal., accepted, DOI:
  10.1007/s12220-011-9275-z.

\bibitem{svdm-tpcurves}
Pawe{\l} Strzelecki and Heiko von~der Mosel.
\newblock Tangent-point self-avoidance energies for curves.
\newblock {\em J. Knot Theory Ramifications}, 21(5):28 pages, 2012.

\bibitem{sullivan_2002}
John~M. Sullivan.
\newblock Approximating ropelength by energy functions.
\newblock In {\em Physical knots: knotting, linking, and folding geometric
  objects in {$\Bbb R^3$} ({L}as {V}egas, {NV}, 2001)}, volume 304 of {\em
  Contemp. Math.}, pages 181--186. Amer. Math. Soc., Providence, RI, 2002.

\bibitem{tolsa-survey}
Xavier Tolsa.
\newblock Analytic capacity, rectifiability, and the {C}auchy integral.
\newblock In {\em International {C}ongress of {M}athematicians. {V}ol. {II}},
  pages 1505--1527. Eur. Math. Soc., Z\"urich, 2006.

\bibitem{verdera_2001}
Joan Verdera.
\newblock {$L^2$} boundedness of the {C}auchy integral and {M}enger curvature.
\newblock In {\em Harmonic analysis and boundary value problems
  ({F}ayetteville, {AR}, 2000)}, volume 277 of {\em Contemp. Math.}, pages
  139--158. Amer. Math. Soc., Providence, RI, 2001.

\end{thebibliography}
\bibliographystyle{plain}

}

\medskip\noindent
\textsc{Pawe\l{} Strzelecki, Marta Szuma\'{n}ska}: Instytut
Matematyki, Uniwersytet Warszawski, ul. Banacha 2, PL-02-097
Warsaw, Poland.

\medskip\noindent
{\sc Heiko von der Mosel}: Institut f\"ur Mathematik, RWTH Aachen
University, Templergraben 55, D-52062 Aachen, Germany.

\end{document}